\documentclass{amsart}
\usepackage{amsmath, amssymb}
\usepackage{url}
\input{xypic}

\newtheorem{theorem}{Theorem}[section]

\newtheorem{lemma}[theorem]{Lemma}
\newtheorem{corollary}[theorem]{Corollary}
\newtheorem{fact}[theorem]{Fact}
\newtheorem{proposition}[theorem]{Proposition}
\newtheorem{claim}[theorem]{Claim}

\theoremstyle{definition}

\newtheorem{question}[theorem]{Question}
\newtheorem{remark}[theorem]{Remark}
\newtheorem{definition}[theorem]{Definition}

\newenvironment{claimproof}[1]{\par \vspace{.1in} \noindent \emph{Proof of Claim:}\space#1}{\hfill $\maltese$ \vspace{.1in}}

\def \U {\mathcal U}
\def \C {\mathcal C}
\def \M {\mathcal M}
\def \K {\mathcal K}
\def \G {\mathcal G}
\def \H {\mathcal H}
\def \A {\mathcal A}
\def \F {\mathcal F}
\def \O {\mathcal O}
\def \I {\mathcal I}

\def\GL{\operatorname{GL}}
\def\Pic{\operatorname{Pic}}
\def\Sp{\operatorname{Sp}}
\def\PGL{\operatorname{PGL}}
\def\dcf{\operatorname{DCF}}
\def\acf{\operatorname{ACF}}
\def\dccm{\operatorname{DCCM}}
\def\ccm{\operatorname{CCM}}

\def\dcl{\operatorname{dcl}}
\def\acl{\operatorname{acl}}

\def\id{\operatorname{id}}

\def\rank{\operatorname{rank}}

\def\tp{\operatorname{tp}}

\def\aut{\operatorname{Aut}}

\def\loc{\operatorname{loc}}
\def\trdeg{\operatorname{trdeg}}

\def \PP {\mathbb P}
\def \AA {\mathbb A}
\def \CC {\mathbb C}

\def\Ind#1#2{#1\setbox0=\hbox{$#1x$}\kern\wd0\hbox to 0pt{\hss$#1\mid$\hss}
\lower.9\ht0\hbox to 0pt{\hss$#1\smile$\hss}\kern\wd0}
\def\ind{\mathop{\mathpalette\Ind{}}}
\def\Notind#1#2{#1\setbox0=\hbox{$#1x$}\kern\wd0\hbox to 0pt{\mathchardef
\nn=12854\hss$#1\nn$\kern1.4\wd0\hss}\hbox to
0pt{\hss$#1\mid$\hss}\lower.9\ht0 \hbox to
0pt{\hss$#1\smile$\hss}\kern\wd0}

\title[Definable Galois theory]{Definable Galois theory\\ for bimeromorphic geometry}
\author{Rahim Moosa}
\address{Rahim Moosa, Department of Pure Mathematics, University of Waterloo, Canada}
\author{Anand Pillay}
\address{Anand Pillay, Department of Mathematics, University of Notre Dame, USA}
\makeatletter

\let\@wraptoccontribs\wraptoccontribs
\makeatother
\contrib[with appendices by]{Frederic Campana and Matei Toma}

\date{\today}
\thanks{The first author was supported by an NSERC Discovery Grant.
The second author was supported by NSF grants DMS-2054271 and DMS-2502292.}
\begin{document}

\begin{abstract}
The outlines of a ``Galois theory" for bimeromorphic geometry is here developed, via the study of model-theoretic definable binding groups in the theory $\ccm$ of compact complex spaces.
As an application, a structure theorem about principal meromorphic bundles with algebraic structure group, and admitting no horizontal subvarieties, is deduced.
Examples of algebraic groups arising as binding groups are provided, as is a characterisation of when they are linear.
Using binding groups in $\ccm$ it is shown that, in contrast to the situation in differentially closed fields, there are many algebraic groups which admit nontrivial definable torsors over $\acl$-closed sets in the theory $\dccm$ of existentially closed differential $\ccm$-structures.
A self-contained exposition of the binding group theorem in totally transcendental theories, that emphasises the bitorsorial nature of the construction, is also included.
\end{abstract}


\maketitle

\setcounter{tocdepth}{1}
\tableofcontents

\vfill\pagebreak

\section{Introduction}

\noindent
The geometric setting for this article is what might be called the {\em meromorphic category} whose objects are of the form $X:=\overline X\setminus E$ where $\overline X$ is a reduced and irreducible compact complex analytic variety and $E\subset\overline X$ is a proper closed analytic subset, and the morphisms are dominant holomorphic maps on~$X$ that extend to meromorphic maps on the compactification~$\overline X$.
Any such~$X$ is endowed with a noetherian topology induced by the closed analytic subsets of~$\overline X$, which we call the {\em Zariski topology} on~$X$.
The meromorphic category is an expansion of quasi-projective complex algebraic geometry, and, by Chow's algebraicity theorem, it induces no additional structure on algebraic varieties.

The meromorphic category is accessible to model theory via the first order multi-sorted structure where there is a sort for each such~$X$ and a predicate for each Zariski closed subset of a finite cartesian product of sorts.
The theory of this structure, denoted by $\ccm$, is tame: it is totally transcendental, of finite Morley rank (sort-by-sort), admits quantifier elimination, and admits elimination of imaginaries.
See~\cite{ccs} for a survey of the basic model theory of $\ccm$, including references for the above facts.
In Section~\ref{sec:ccmprelims} below we give a quick recap and overview of the aspects that will concern us.

We denote by~$\C$ the field definable in $\ccm$ whose interpretation in the standard model is the complex field; it lives on the sort of the complex affine line.
The induced structure on~$\C$ is that of an algebraically closed field with no additional structure.
This is the model-theoretic interpretation of quasi-projective algebraic geometry as a reduct of meromorphic geometry.
Our work here has to do with using model theory to analyse parts of the full meromorphic category in terms of this algebro-geometric reduct.
More specifically, we work out in $\ccm$ various aspects of the general theory of {\em binding groups} associated to {\em $\C$-internality}.
A typical geometric situation exhibiting these model-theoretic  phenomena is  a surjective holomorphic map of compact complex analytic varieties, $f:X\to Y$,  such that after base extension~$X$ embeds meromorphically into a product of (the new) base with some complex projective space.
The binding group in this case is an algebraic family of algebraic groups over~$Y$ acting meromorphically on the fibres of~$X$ over~$Y$.
If we assume further that $f^*:\CC(Y)\to\CC(X)$ is an isomorphism (that is there are no new meromorphic functions on~$X$), then this group action is transitive on the fibres.
If we assume that $Y$ admits no nonconstant meromorphic functions (that is the algebraic dimension of~$Y$ is zero) then the binding group is a fixed complex algebraic group acting meromorphically on~$X$ over~$Y$.

While $\ccm$ binding groups have been mentioned in a few places (for example, \cite[$\S$5]{c3}, \cite[$\S$8]{nmdeg}, \cite[$\S$4.3]{jjp}, \cite[$\S$7]{abred}), the subject remains largely unexplored, compared, for example, to binding groups in the theory of differentially closed fields ($\dcf_0$), which has a number of striking similarities with $\ccm$.
In $\dcf_0$, binding groups (relative to the field of constants) give a model-theoretic account (and generalisation) of the differential Galois theories of Picard-Vessiot and Kolchin.
We are not aware of a ``classical" Galois theory in the meromorphic category; indeed, we see part of our work here as suggesting the outlines of such a theory.

Our main contributions are as follows:

\medskip

(1)
We take the opportunity, in Section~\ref{sec:bgr} below, to include yet another exposition of the general model theory of binding groups in totally transcendental theories.
Binding groups (or definable groups of automorphisms) first make an appearance in Zilber's work on the fine structure of uncountably categorical theories.
See~\cite{Zilber} for an account,  as well as references to his  original papers.
Poizat~\cite{Poizat} gives a model-theoretic account of Kolchin's differential Galois theory, via binding groups.
Hrushovski, in his Ph. D. thesis, gives an account in the context of stable theories (with important applications).
Expositions are also  is given in~\cite{stablegroups} and~\cite{gst}.
The theory is eventually put into a very general and comprehensive form by Hrushovski in~\cite{hrushovski-bgt}. 
Our exposition, which remains restricted to totally transcendental  theories, emphasises the ``bitorsorial" nature and canonicity of the construction, with (we believe) somewhat simpler proofs.
We build on and simplify the treatment in~\cite{Leon-Sanchez-Pillay}. 

\medskip

(2)
As an application of the theory of binding groups in~$\ccm$ we deduce the following structural result in meromorphic geometry, characterising when principal bundles for algebraic groups admit no ``horizontal" subvarieties:

\begin{theorem}
\label{1}
Suppose $G$ is a (not necessarily connected) complex algebraic group and $f:X\to Y$ is a principal meromorphic $G$-bundle, with $Y$ of algebraic dimension~$0$.
Then the following are equivalent:
\begin{itemize}
\item[(i)]
$X$ admits no proper Zariski closed subsets projecting dominantly onto~$Y$.
\item[(ii)]
$f$ admits no principal meromorphic $H$-subbundles, for~$H$ a proper complex algebraic subgroup of~$G$.
\end{itemize}
\end{theorem}

This appears as Corollary~\ref{cor:tt}, below, and we point the reader to Sections~\ref{sec:ccmprelims} and~\ref{sec:bundle} for a precise explanation of the terminology.
 What the assumptions, that $f:X\to Y$ is a principal meromorphic $G$-bundle and~$Y$ is of algebraic dimension~$0$, give us is a complex algebraic subgroup~$G'\leq G$ that is an avatar of the binding group of the generic fibre of~$f$.
 Conditions~(i) and~(ii) are then shown to both be equivalent to $G'=G$.
The general strategy of proof here works even if we drop the assumption that the algebraic dimension of~$Y$ is zero, however in that case~$G'$ will not necessarily be defined over~$\CC$, it will be defined over (a finite extension of) the field $\CC(Y)$ of meromorphic functions on~$Y$.
An articulation of the theorem in that general case would thus involve a version of condition~(ii) that replaces the complex algebraic subgroup~$H$ by a moving family of algebraic subgroups of~$G$ over finite covers of~$Y$, and is not presented here.
%A version of Theorem~\ref{1} could therefore be articulated in the more general case, but it would require making sense of principle subbundles with respect to such~$H$.

Theorem~\ref{1} bears some resemblance to an old theorem of Campana's, namely~\cite[Corollaire~3]{campana}.
The precise relationship remains somewhat unclear.
Under the assumption that~$X$ is of K\"ahler-type, we can deduce Theorem~\ref{1} from Campana's theorem, though still using some model theory (isolated types and prime models).
Conversely, some aspects of Campana's theorem are recovered by Theorem~\ref{1}.
It would be of real interest to compare the binding group action we develop here with the relative meromorphic group action constructed in~\cite[Corollaire~3]{campana},  but we delay this to future work.

\medskip

(3)
We produce a few new examples of algebraic groups that arise as binding groups in $\ccm$; namely,  simple abelian varieties, the multiplicative torus, and certain algebraic subgroups of $\PGL_{2n}$ of arbitrarily high dimension.
(The only specific previously presented example was $\PGL_2$ in~\cite{jjp}.)
These appear in Sections~\ref{subsect:abelian} and~\ref{subsect:gm}, and are based on material in the two included appendices by Fr\'ed\'eric Campana and Matei Toma, respectively.

\medskip

(4)
We characterise when binding groups are linear:

\begin{theorem}
Suppose~$p$ is a complete type in $\ccm$ that is fundamental $\C$-internal and weakly $\C$-orthogonal.
The binding group of~$p$ is a linear algebraic group  if and only if~$p$ is interdefinable with the type of a basis for a definable $\C$-vector space.
\end{theorem}

This appears as Propositions~\ref{prop:getvs} and~\ref{vs-linalg} below.
We point the reader to Section~\ref{sec:bgr} for precise definitions of fundamental internality and weak orthogonality, and to Section~\ref{subsect:intorth} for their manifestations in the meromorphic category.
The conclusion here is not all that one could desire: a concrete way for definable $\C$-vector spaces to arise is from relative linear spaces in the standard model, and we ask in Question~\ref{question} below if these are the only ones.

When~$p$ is the generic type of the generic fibre of a holomorphic map $f:X\to Y$ between K\"ahler-type spaces, we also obtain a {\em geometric} characterisation of when the binding group is linear; namely that the relative Albanese of~$f$ be trivial.
This is done in Proposition~\ref{prop:linalgbg} below, to which we refer the reader for a precise statement.

\medskip

(5)
We use binding groups in $\ccm$ to show that a certain natural, and possibly expected extension of a result of Kolchin's, does {\em not} hold.
The original theorem we have in mind is from~\cite{dag} and says that over an algebraically closed differential field every $\dcf_0$-definable torsor for an algebraic group is trivial.
This plays a significant role in differential Galois theory.
See~\cite{dgt1} for generalisations from a model-theoretic perspective.
Now, a theory $\dccm$ was introduced in~\cite{dccm} that expands $\ccm$ in precisely the same way that $\dcf_0$ expands $\acf_0$.
As such it is natural to ask whether, over an $\acl$-closed set of parameters, every $\dccm$-definable torsor of a $\ccm$-definable group is trivial.
We will show that this is false.
Indeed, every positive-dimensional algebraic group arising as a binding group in $\ccm$, or even containing a binding group in $\ccm$, admits a nontrivial $\dccm$-definable torsor.
This is carried out in Section~$\S$\ref{sect:h1}.

\medskip

\subsection*{Acknowledgements}
The second author would like to thank David Meretzky for helping him understand the bitorsorial nature of the binding group construction, Gareth Jones for asking some
pertinent questions (answered by our Theorem~\ref{1}), and Remi Jaoui for catching an error in an earlier version of this paper.
Both authors are indebted to Fr\'ed\'eric Campana and Matei Toma, for useful discussions, and for agreeing to contribute appendices to this paper.
Finally, the authors are grateful for the hospitality of the Fields Institute in Toronto where some of this work was done.

\bigskip
\section{Binding groups revisited}
\label{sec:bgr}

\noindent
In this section we give a summary of the binding group construction in the totally transcendental case.
The purpose is mostly to fix notation and terminology, but we also take the opportunity to make explicit the bitorsorial nature of the construction.
The latter is well known in the differential context (see~\cite{bertrand}) and also appears, somewhat differently and in greater generality, in~\cite[Theorem~B.1]{hrushovski-bgt}.

We work in a sufficiently saturated model $\overline\M$ of a complete (multi-sorted) totally transcendental theory~$T$ which eliminates imaginaries and admits quantifier elimination.
All our parameter sets are assumed to be {\em small}, namely of cardinality strictly less than the level of saturation, unless explicitly stated otherwise.
We fix a base set of parameters~$B$, a partial type $\Phi$ over~$B$, and a $B$-definable set~$\C$.
We will be focusing eventually on the case when $\Phi$ is a complete type $p\in S(B)$, but for now we work generally.
The central notion is that of internality:

\begin{definition}
We say that~$\Phi$ is {\em $\C$-internal} if $\Phi(\overline\M)\subseteq\dcl(B'\C)$ for some (small) $B'\supseteq B$.
\end{definition}

We are interested in the (abstract) group of elementary permutations of
$\Phi(\overline\M)$ over $B\cup\C$, as well as its definable avatars (which exist under the assumption of $\C$-internality).

\begin{definition}
We denote by $\aut_B(\Phi/\C)$ the group of restrictions to $\Phi(\overline\M)$ of elements of $\aut_{B\C}(\overline\M)$.
By a {\em definable binding group action for~$\Phi$ relative to~$\C$ over~$B$}, sometimes described simply as a {\em binding group}, we mean a $B$-definable group, $\G$, together with a  (relatively) $B$-definable faithful left action on $\Phi(\overline\M)$, such that there is an (abstract) group isomorphism between~$\G$ and~$\aut_B(\Phi/\C)$ respecting the action on $\Phi(\overline\M)$.
\end{definition}

\begin{remark}
\label{rem:canonicaliso}
The isomorphism between $\aut_B(\Phi/\C)$ and a definable binding group~$\G$ is canonical; it must take $\sigma\in\aut_B(\Phi/\C)$ to the unique element $g\in\G$ that has the same action on $\Phi(\M)$ as $\sigma$.
\end{remark}

The binding group theorem -- as it appears, for example, in the very general~\cite[Theorem~B.1']{hrushovski-bgt}, but then specialised to totally transcendental theories -- asserts that if~$\Phi$ is $\C$-internal then a definable binding group action for~$\Phi$ relative to~$\C$ exists.
(When $\Phi$ is complete and stationary, see also~\cite[$\S$7.4]{gst}.)

In the case of complete types the following conditions are of central importance:

\begin{definition}
Suppose~$p\in S(B)$.
We say that~$p$ is {\em fundamental $\C$-internal} if $p(\overline\M)\subseteq\dcl(Ba\C)$ for some (equivalently any) $a\models p$.
We say that~$p$ is {\em weakly $\C$-orthogonal} if it isolates a complete type over $B\cup\C$.
\end{definition}

\begin{remark}
Weak $\C$-orthogonality implies, in particular, that if $a\models p$ and~$c$ is a tuple from~$\C$ then $a\ind_Bc$.
When~$p$ is stationary, the latter condition characterises weak $\C$-orthogonality and is often used as its definition.
\end{remark}

Note that fundamental $\C$-internality ensures the action of $\aut_B(p/\C)$ on $p(\overline\M)$ is free, while weak $\C$-orthogonality makes the action transitive.
Another consequence of weak $\C$-orthogonality (in the presence of $\C$-internality) is that $p$ is isolated; see, for example, \cite[Lemma~2.2]{Leon-Sanchez-Pillay}.
The upshot is that if~$p\in S(B)$ is fundamental $\C$-internal and weakly $\C$-orthogonal, and $\G$ is a definable binding group action for~$p$ relative to~$\C$ over~$B$, then $p(\overline\M)$ is a {\em $B$-definable torsor} for~$\G$, namely it is a $B$-definable set admitting a $B$-definable uniquely transitive action of~$\G$.

\begin{remark}
\label{rem:uniquebg}
In the fundamental and weakly orthogonal case, binding group actions are uniquely determined up to $B$-definable isomorphism.
Indeed, any two $B$-definable uniquely transitive group actions on the same set, say $\G$ and $\G'$ on $P$, are canonically $B$-definably isomorphic by identifying $g\in \G$ with $g'\in \G'$ if $ga=g'a$ for some (equivalently any) $a\in P$.
We therefore talk about {\em the} binding group in this case.
\end{remark}

Recall the following useful characterisation of weak orthogonality:

\begin{lemma}
\label{lem:wochar}
$\tp(a/B)$ is weakly $\C$-orthogonal if and only if
$$\dcl(Ba)\cap\dcl(B\C)=\dcl(B).$$
\end{lemma}

\begin{proof}
($\Longleftarrow$).
Stable embeddedness of $\C$ implies that $\tp(a/B\C)$ is isolated by $\tp(a/BE)$ where $E:=\dcl(Ba)\cap\dcl(\C)$.
(See~\cite[Appendix]{acfa} for a study of stable embeddedness.)
Our assumption implies that $E\subseteq\dcl(B)$, so that $\tp(a/B)$ isolates $\tp(a/B\C)$, as desired.

($\Longrightarrow$).
Suppose $d\in\dcl(B\C)$ is of the form $d=f(a)$ for some $B$-definable function $f$.
Let $\sigma\in\aut_B(\overline\M)$.
Then $\sigma(d)=f(\sigma(a))$.
On the other hand, our assumption implies that~$a$ and $\sigma(a)$ have the same type over $B\cup\C$, and hence over $\dcl(B\C)$.
So $f(\sigma(a))=d$ as well.
That is, $\sigma(d)=d$.
Hence $d\in\dcl(B)$, as desired.
\end{proof}

The following to a large extent justifies our focus on the case of complete types that are fundamental $\C$-internal and weakly $\C$-orthogonal.

\begin{proposition}
\label{prop:redtofundwo}
Suppose $\Phi$ is a partial type over~$B$ that is $\C$-internal.
There exists an $n$-tuple of realisations $\bar a$ of $\Phi$, and a finite tuple $c\subset\dcl(B\bar a)\cap\dcl(\C)$, such that, setting $q:=\tp(\bar a/Bc)$,
\begin{itemize}
\item[(a)]
$q$ is fundamental $\C$-internal,
\item[(b)]
$q$ is weakly $\C$-orthogonal,
\item[(c)]
there is an isomorphism $\alpha:\aut_B(\Phi/\C)\to\aut_{Bc}(q/\C)$ given by $\alpha(\sigma)=\beta$ if and only if~$\beta$ agrees with the diagonal action of~$\sigma$ on $q(\overline\M)\subseteq p(\overline\M)^n$,
\item[(d)]
$\alpha$  induces a $Bc$-definable isomorphism between any binding group for $\Phi$ relative to~$\C$ over~$B$ and the binding group for~$q$ relative to~$\C$ over $Bc$.
\end{itemize}
\end{proposition}

\begin{proof}
By compactness we can find a finite tuple~$b$ such that $\Phi(\overline\M)\subseteq\dcl(Bb\C)$.
We first argue that we can replace~$b$ by a tuple of realisations of~$\Phi$.
Indeed, consider $\tp(b/B\C\Phi(\overline\M))$.
Using the stable embeddedness of $\Phi(\overline\M)$, if we let
$$E:=\dcl(B\C b)\cap \dcl(\Phi(\overline\M))$$
then $\tp(b/B\C\Phi(\overline\M))$ is isolated by $\tp(b/B\C E)$.
By total transcendentality,
$$\dcl(B\C E)\subseteq\dcl(B\C\bar a)$$
for some finite tuple $\bar a\subset\Phi(\overline\M)$.
Hence $\tp(b/B\C\bar a)$ isolates $\tp(b/B\C\Phi(\overline\M))$.
Since $\Phi(\overline\M)\subseteq\dcl(Bb\C)$, it follows that $\Phi(\overline\M)\subseteq\dcl(B\C\bar a)$ also.

Next, consider $\tp(\bar a/B\C)$ and use the stable embeddedness of~$\C$ to obtain, exactly as above, finite~$c\subset\dcl(B\bar a)\cap\dcl(\C)$ such that 
 $q:=\tp(\bar a/Bc)$ isolates $\tp(\bar a/B\C)$.
In particular $q$ is weakly $\C$-orthogonal.
But also, $q(\overline\M)\subseteq \dcl(\Phi(\overline\M))\subseteq\dcl(B\C \bar a)$, witnessing that $q$ is fundamental $\C$-internal.

Restriction induces surjective homomorphisms from $\aut_{B\C}(\overline\M)$ to $\aut_B(\Phi/\C)$ and $\aut_{Bc}(q/\C)$, the latter using that the additional parameters~$c$ are from~$\dcl(\C)$.
Since $\bar a$ is a tuple of realisations of~$\Phi$, we have that $Q:=q(\overline\M)$ is contained in a cartesian power of $P:=\Phi(\overline\M)$, so that these restrictions induce a surjective homomorphism from $\aut_B(\Phi/\C)$ to $\aut_{Bc}(q/\C)$,
$$\xymatrix{
& \aut_{B\C}(\overline\M)\ar[dr]\ar[dl] &\\
\aut_B(\Phi/\C)\ar@{-->}[rr]^\alpha & & \aut_{Bc}(q/\C)
}$$
which is in fact an isomorphism as $P\subseteq\dcl(B\C\bar a)$.
Given definable binding groups for~$\Phi$ and~$q$ relative to~$\C$, say~$\G$ and~$\H$ respectively, we obtain an induced isomorphism $\widehat\alpha:\G\to\H$, 
where $\widehat\alpha(g)=h$ if and only if~$h$ agrees with the diagonal action of~$g$ on $Q$.
This isomorphism is $Bc$-definable as~$q$ is isolated, and hence~$Q$ is a $Bc$-definable set.
\end{proof}

\begin{remark}
In Proposition~\ref{prop:redtofundwo}, even if we began with~$\Phi$ a complete and stationary type, it is not necessarily the case that~$q$ will be stationary.
Indeed, $q$ will be stationary precisely if the binding group is connected.
\end{remark}

From now on, we therefore focus on complete types $\Phi=p\in S(B)$ that are fundamental $\C$-internal and weakly $\C$-orthogonal.

There is another definable avatar of $\aut_B(p/\C)$, but this one acting on the right, that has the added benefit of living in the induced structure on~$\C$, and that will be of particular interest to us.
While it is possible to deduce its existence {\em a posteriori} from the binding group theorem, we will instead do the construction the other way round, in the proof of Theorem~\ref{bgt} below.

\begin{definition}
A {\em bitorsor} is a triple $(G,S,H)$ where
\begin{itemize}
\item[(i)]
$G$ is a group with a uniquely transitive left action on the set~$S$,
\item[(ii)]
$H$ is a group with a uniquely transitive right action on~$S$, such that
\item[(iii)]
$(ga)h=g(ah)$ for all $g\in G, a\in S, h\in H$.
\end{itemize}
\end{definition}

\begin{remark}
\label{rem:bitorsor}
If $(G,S,H)$ is a bitorsor then~$H^{\operatorname{opp}}$ naturally acts on~$S$ on the left, and is precisely the group of permutations of~$S$ that commute with the action of~$G$.
Note also that each $a\in S$ gives rise to a group isomorphism $\rho_a:G\to H$ given by $ga=a\rho_a(g)$ for each $g\in G$.
\end{remark}

\begin{theorem}[Binding group theorem, bitorsor version]
\label{bgt}
Suppose $p\in S(B)$ is fundamental $\C$-internal and weakly $\C$-orthogonal.
Then there exists a $B$-definable bitorsor $(\G, p(\overline\M),\widehat\G)$ such that
\begin{itemize}
\item[(a)] $\G$ is the binding group for~$p$ relative to~$\C$ over~$B$, and
\item[(b)] $\widehat\G\subseteq\dcl(B\C)$.
\end{itemize}
\end{theorem}

\begin{proof}
We first construct $\widehat\G$ directly and use it to construct $\G$.

Let $P:=p(\overline\M)$.
As~$p$ is isolated, $P$ is $B$-definable.
And, the action of $\aut_B(p/\C)$ on $P$ is uniquely transitive.
These facts will be used repeatedly in the arguments that follow.

By compactness (or rather saturation), fundamental internality is witnessed by a $B$-definable partial function, $f(-,-)$, such that for any $a,a'\in P$, there is a tuple $c$ in $\C$ (of fixed appropriate length) such that $a' = f(a,c)$.
Fixing $a\in P$ we define an equivalence relation~$E$ on the appropriate cartesian power of $\C$ by setting
$c E c'$ if  $f(a,c) = f(a,c')$.
As any two elements of $P$ have the same type over $B\cup\C$, $E$ does not depend on the choice of~$a$, and is $B$-definable.
Without loss (replacing $c$ by $c/E$) we may assume that for each $a,a'\in P$ there is a unique $c\in\dcl(B\C)$ such that $a' = f(a,c)$.
The set of  $c$ that arise in this way is then a $B$-definable set $X$ living in $\dcl(B\C)$.

For any $a\in P$ and $c\in X$, let us write  $a\cdot c\in P$ for $f(a,c)$.
Fixing $a\in P$, we have that for any $c_{1}, c_{2}\in X$ there is unique $c_{3}\in X$ such that  $(a\cdot c_{1})\cdot c_{2} = a\cdot c_{3}$.
Again, the function that assigns $c_3$ to $(c_1,c_2)$ does not depend on~$a$ as any two elements of~$P$ have the same type over $B\cup\C$.
We have, thus, a $B$-definable binary operation on $X$, that we will also denote by $\cdot$, having the property that $(a\cdot c_1)\cdot c_2 = a\cdot(c_1\cdot c_2)$ for any $a\in P$ and $c_1,c_2\in X$.

\begin{claim}
\label{claim:h}
Fix $a_\circ\in P$.
For  $\sigma\in\aut_B(p/\C)$,  let $c_{\sigma}$ denote the unique element of~$X$ such that
$\sigma(a_\circ) = a_\circ\cdot c_{\sigma}$.
Then $\sigma\mapsto c_{\sigma}$ is a bijection between $\aut_B(p/\C)$ and~$X$, satisfying $ c_{\sigma\circ\tau} =  c_{\sigma}\cdot c_{\tau}$.
\end{claim} 

\begin{claimproof}
Injectivity:
If $c_\sigma=c_\tau$ then $\sigma(a_\circ)=a_\circ\cdot c_\sigma=a_\circ\cdot c_\tau=\tau(a_\circ)$, and hence $\sigma=\tau$ by the freeness of the action.
Surjectivity: given $c\in X$, transitivity of the action gives $\sigma\in\aut_B(p/\C)$ such that $\sigma(a_\circ)=a_\circ\cdot c$, and so $c=c_\sigma$.

Finally, we compute:
\begin{equation*}
\label{comp1}
(\sigma\circ \tau)(a_\circ) = \sigma(\tau(a_\circ)) = \sigma(a_\circ\cdot c_{\tau}) = \sigma(a_\circ)\cdot c_{\tau}=(a_\circ\cdot c_{\sigma})\cdot c_{\tau} = a_\circ\cdot(c_{\sigma}\cdot c_{\tau}),
\end{equation*}
where the third equality is because $c_\tau$ is from $\C$ and hence fixed by $\sigma$, and the last equality is by definition of $\cdot$ on $X$.
\end{claimproof}

It follows from Claim~\ref{claim:h} that $\widehat\G:=(X,\cdot)$ is a $B$-definable group and $\sigma\mapsto c_{\sigma}$ is an isomorphism between $\aut_B(p/\C)$ and $\widehat\G$. (This isomorphism does depend, however, on some choice of $a_\circ$.)
Working in the group $\widehat\G$ we now drop the symbol $\cdot$ and simply write the group operation as concatenation.

The right action of $\widehat\G$ on $P$ is given by $(a,c)\mapsto a\cdot c=f(a,c)$, for all $a\in P$ and $c\in X$.
It is clearly a $B$-definable function.
That it is a right group action is exactly the defining property of the group operation of~$\widehat\G$.
That it is a uniquely transitive action follows from the fact that for each $a,a'\in P$ there is a unique $c\in X$ such that $a' = f(a,c)$.

At this point it is worth noting that we already have a bitorsor $(\aut_B(p/\C), P,\widehat\G)$.
Indeed, as $X\subseteq\dcl(B\C)$ it is fixed pointwise by elements of $\aut_B(p/\C)$.
Hence $\sigma(a)\cdot c=\sigma(a\cdot c)$ for all $\sigma\in\aut_B(p/\C)$, $a\in P$, and $c\in X$, which is precisely the bitorsor condition.

What remains is to find a $B$-definable avatar of the left hand side, namely the torsor 
$(\aut_B(p/\C), P)$.
By the freeness of the action of $\aut_B(p/\C)$, to every pair $(a,a')\in P\times P$ there is a unique $\sigma\in\aut_B(p/\C)$ such that $\sigma(a)=a'$.
So we have a well-defined surjective function $P\times P\to\aut_B(p/\C)$.
Let $E$ be the equivalence relation on $P\times P$ given by the fibres of this function.
We claim that $E$ is $B$-definable.
Indeed, this follows once we observe that $(a,a')E(b,b')$ if and only if there is (unique) $c\in X$ such that $a\cdot c=b$ and $a'\cdot c=b'$.
That left implies right is because $(a,a')E(b,b')$ implies $\tp(ab/B\C)=\tp(a'b'/B\C)$, and so if we choose~$c\in X$ such that $a\cdot c=b$, then as $c\in\dcl(B\C)$ it must also be that $a'\cdot c=b'$.
For the converse, choosing $\sigma\in\aut_B(p/\C)$ such that $\sigma(a)=a'$ we have that 
$$\sigma(b)=\sigma(a\cdot c)=\sigma(a)\cdot c=a'\cdot c=b'$$
so that $(a,a')E(b,b')$, as desired.

We thus have a $B$-definable set $Y:=(P\times P)/E$ which admits a bijection to $\aut_B(p/\C)$ given by sending the $E$-class of $(a,a')$ to the unique $\sigma\in\aut_B(p/\C)$ such that $\sigma(a)=a'$.
This bijection endows~$Y$ with a definable group structure, that we denote by~$\G$, and with a left action of~$\G$ on~$P$.
It only remains to show that the group and its action are $B$-definable.
But the action of~$\G$ on~$P$ satisfies $y\cdot a=a'$ if and only if $(a,a')\in y$, which is visibly $B$-definable.
Since the group structure can be recovered from the faithful action (namely, $y_1y_2 = y_3$ if and only if $y_1\cdot (y_2\cdot a) = y_3\cdot a$ for all $a\in P$) this too is $B$-definable.
\end{proof} 

Just as the binding group action is uniquely determined upto $B$-definable isomorphism (see Remark~\ref{rem:uniquebg}) so too is the group $\widehat\G$ given by Theorem~\ref{bgt}, along with its right action on $p(\overline\M)$.
This allows us to name it:

\begin{definition}
Suppose $p\in S(B)$ is fundamental $\C$-internal and weakly $\C$-orthogonal.
We call the group $\widehat\G$ given by Theorem~\ref{bgt} the {\em dual binding group of~$p$ relative to $\C$}.
That is, $\widehat\G$ is a $B$-definable group living in $\dcl(B\C)$ with a $B$-definable right action on $p(\overline\M)$ that forms, together with the left action of the binding group~$\G$,  a bitorsor $(\G,p(\overline\M),\widehat\G)$.
\end{definition}

By Remark~\ref{rem:bitorsor}, $\widehat\G^{\operatorname{opp}}$, together with its induced left-action on $p(\overline\M)$, is precisely the group of permutations of $p(\overline\M)$ that commutes with the action of the binding group~$\G$.
In other words, $\widehat\G^{\operatorname{opp}}$ is what Hrushovski calls the ``opposite" binding group and denotes by $\aut_{\G}(p(\overline\M))$ in his treatment of the binding group theorem in~\cite[Theorem~B.1]{hrushovski-bgt}.

We note that Theorem~\ref{bgt} has a converse:

\begin{remark}
Suppose $(G,S,H)$ is a bitorsor defined over~$B$ such that
\begin{itemize}
\item[(i)]
The action of~$G$ on~$S$ preserves types over $B\cup{\mathcal C}$, and 
\item[(ii)]
$H\subseteq\dcl(B\C)$.
\end{itemize}
Then $S=p(\M)$ for some $p\in S(B)$ that is fundamental ${\mathcal C}$-internal and weakly $\C$-orthogonal, $G$ is the binding group for~$p$ relative to~$\C$ over~$B$, and $H$ is the dual binding group. 
\end{remark}

\begin{proof}
As the action of~$G$ on~$S$ preserves types over $B\cup{\mathcal C}$ it follows that all elements of $S$ have the same type over $B\cup{\mathcal C}$, which, as~$S$ is $B$-definable, implies that~$S=p(\M)$ for some $p\in S(B)$ which is weakly orthogonal to~${\mathcal C}$. 
On the other hand, as~$S$ is a definable torsor for~$H$, (ii) implies that~$p$ is fundamental ${\mathcal C}$-internal. 

To see that~$G$ is the~binding group of $p$ (over $B$ and relative to ${\mathcal C}$) we should check that every $\sigma\in\aut_{B}(p/{\mathcal C})$ is given by some $g\in G$.
Let $\sigma\in G$ and fix $a\in S$.
Since $\sigma(a)\in S$ we have that $\sigma(a) = ga$ for some $g\in G$. We claim that $\sigma(b) = gb$ for all $b\in S$.
Indeed, given $b$ and letting $h\in H$ be such that  $b = ah$, we have that
$\sigma(b) = \sigma(ah) = \sigma(a)h = (ga)h = g(ah) = gb$, as desired.

That~$H$ is the dual binding group now follows from the definition.
\end{proof}

Galois groups in field theory, and in differential algebra, arise as automorphism groups of (differential) field extensions generated by solutions to a (differential) polynomial equation.
%In field theory, the Galois group of a polynomial equation refers, roughly speaking, to the group of automorphisms of the field extension generated by a solution.Similarly, one has differential Galois groups associated to certain algebraic differential equations in differential algebra.
We now point out how, in the general setting of totally transcendental theories, the binding group and its dual, as given by Theorem~\ref{bgt}, can also be viewed this way.
%relate to the abstract analogue of such Galois groups.
The role of the (differential) field extension is played by an extension of definably closed substructures $B\subseteq D$ where $D=\dcl(Ba)$ and~$a$ is a realisation of the type in question.
We use $\aut_B(D)$ to denote the group of automorphism of the structure~$D$ that fix~$B$ pointwise.
Also, if~$X$ is a $D$-definable set then we denote by $X(D)$ the points of~$X$ which are in~$D$.

\begin{corollary}
\label{cor:extrinsic}
Suppose $p\in S(B)$ is weakly $\C$-orthogonal and fundamental $\C$-internal where $B$ is definably closed.
Let $\G$ be the binding group of $p$ relative to~$\C$ over~$B$, and $\widehat\G$ the dual binding group.
Set $P:=p(\overline\M)$, fix $a\in P$, and consider $D:=\dcl(Ba)$.
\begin{itemize}
\item[(a)]
The action of $\G(D)$ on $P(D)$ is induced by an isomorphism of $\G(D)$ with $\aut_B(D)$.
That is, there is a group isomorphism
$$t:\aut_B(D)\to \G(D)$$
such that $\sigma$ and $t(\sigma)$ act the same on $P(D)$, for all $\sigma\in\aut_B(D)$.
 \item[(b)]
$\widehat\G(D)=\widehat\G(B)$ and
there is a group isomorphism
$$s:\aut_B(D)\to \widehat\G(B)$$
such that $\sigma(a)=as(\sigma)$ for all $\sigma\in\aut_B(D)$.
\end{itemize}
\end{corollary}

\begin{proof}
(a).
Note, first of all, that since~$\G$ acts $B$-definably on~$P$ the restriction is an action of $\G(D)$ on $P(D)$.
The unique transitivity of the action is also inherited.

There is a natural embedding of $\aut_B(D)$ into $\aut_B(p/\C)$ preserving the action on~$P(D)$.
Indeed, by assumption (weak $\C$-orthogonality) we know that $p=\tp(a/B)$ determines a complete type over $B\cup\C$.
Hence so does $\tp(D/B)$.
It follows that every element in $\aut_B(D)$ extends to an element of $\aut_{B\C}(\overline\M)$ which we can then restrict to~$p(\overline\M)$ to get an element of $\aut_B(p/\C)$.
That~$p$ is fundamental $\C$-internal ensures this is well-defined.
It is then clearly an injective group homomorphism.

We have the canonical isomorphism $\aut_B(p/\C)\to\G$ coming from the fact that~$\G$ is a definable binding group (see Remark~\ref{rem:canonicaliso}).
Namely, to $\sigma\in\aut_B(p/\C)$ we associate the unique element $g_\sigma$ such that~$\sigma$ and $g_\sigma$ agree on~$P$.
Let $t$ be the restriction of this to $\aut_B(D)$.
To see that it does map to $\G(D)$, note that if $\sigma\in\aut_B(D)$ then $t(\sigma)=g_\sigma$ is the unique element of $\G$ satisfying $g_\sigma a=\sigma(a)$, and hence $g_\sigma\in\dcl(Ba\sigma(a))\subseteq D$.

It remains only to show that the image of~$t$ is all of $\G(D)$.
To that end, fix $g\in\G(D)$.
Then $g=g_\sigma$ for some $\sigma\in\aut_{B\C}(\M)$.
So $\sigma(a)=ga\in D$.
Hence $\sigma(D)\subseteq D$.
But also $g^{-1}=g_{\sigma^{-1}}$, so that 
$\sigma^{-1}(a)=g^{-1}a\in D$, and hence $\sigma^{-1}(D)\subseteq D$, which implies that $D\subseteq\sigma(D)$.
It follows that $\sigma|_D\in\aut_B(D)$ and $t(\sigma|_D)=g$.

(b).
As $p$ is weakly $\C$-orthogonal, and $B$ is definably closed, Lemma~\ref{lem:wochar} gives us that $\dcl(B\C)\cap D=B$.
Since $\widehat\G\subseteq \dcl(B\C)$, it follows that $\widehat\G(B)=\widehat\G(D)$.

Now, recall the isomorphism $\rho_a:\G\to\widehat\G$ from Remark~\ref{rem:bitorsor}, given by the fact that $(\G,P,\widehat\G)$ is a bitorsor and $a\in P$.
That is, $\rho_a(g)\in\widehat G$ is the unique element such that $ga=a\rho_a(g)$.
This restricts to an isomorphism $\G(D)\to\widehat\G(D)=\widehat\G(B)$.
Precomposing with $t:\aut_B(D)\to \G(D)$ from~(i) yeilds our desired isomorphism $s:\aut_B(D)\to \widehat\G(B)$.
\end{proof}

In differential Galois theory one is used to working in the context where the constants of the base differential field is algebraically closed, and hence any differential closure will admit no new constants.
This assumption also makes sense in the abstract setting, where differential closure is replaced by a prime model.
We record here for future use some consequences: 

\begin{lemma}
\label{fact:widehatbtod}
Suppose $p\in S(B)$ is weakly $\C$-orthogonal and fundamental $\C$-internal,
where $B$ is definably closed,
and let~$\G$ and~$\widehat\G$ be the binding and dual binding groups, respectively.
Fix a prime model $M_B$ over~$B$, let $a\models p$ be a realisation in~$M_B$, and set $P=p(\overline\M)$ and $D=\dcl(Ba)$ as before.
Suppose $\C(M_B)=\C(B)$.
Then
\begin{itemize}
\item[(a)]
$P(M_B) = P(D)$, $\G(M_B)=\G(D)$, and  $\widehat\G(M_B)=\widehat\G(B)$; and
\item[(b)]
restriction induces a surjective homomorphism $\aut_B(M_B)\to\aut_B(D)$.
\end{itemize}
\end{lemma}

\begin{proof}
Note that, as $p$ is isolated there is a realisation~$a$ of~$p$ in~$M_B$.

(a).
As $\widehat\G\subseteq\dcl(B\C)$ and $M_B$ is a model, we have that
$\widehat\G(M_B)\subseteq\dcl(B\C(M_B))$.
From this the third claim follows as $\C(M_B)=\C(B)$.
For the first claim, note that if $a'\in P(M_B)$ then there is $g\in\widehat\G(M_B)=\widehat G(B)$ such that $a'=ga\in P(D)$.
Finally, if $g\in\G(M_B)$ then $g\in\dcl(B,a, ga)$, but $ga\in P(M_B)=P(D)$, so that $g\in \G(D)$.

(b).
To see that restriction to $D=\dcl(Ba)$ makes sense we need to show that if  $a'$ is another realisation of~$p$ in~$M_B$,
then~$a$ and $a'$ are interdefinable over~$B$.
But we have $a'=ag$ for some $g\in\widehat\G(M_B)=\widehat\G(B)$, which witnesses the interdefinability.
For surjectivity, if $\tau\in\aut_B(D)$ then $\tp(a/B)=\tp(\tau(a)/B)$ by quantifier elimination, and hence, as $M_B$ is prime over $B$, there is $\sigma\in\aut_B(M_B)$ such that $\sigma(a)=\tau(a)$, so that $\sigma$ restricts to $\tau$.
\end{proof}

\bigskip
\section{Preliminaries on the theory $\ccm$}
\label{sec:ccmprelims}

\noindent
The usual model-theoretic approach to bimeromorphic geometry (as exposed, for example, in~\cite{ccs}) is to consider the (first order theory of the) structure that has a sort for every reduced and irreducible compact complex analytic space and a predicate for every closed analytic (i.e., {\em Zariski closed}) subset of each finite cartesian power of sorts.
We will, however, follow the slightly more flexible set up of~\cite[$\S$2]{dccm}.

By a {\em meromorphic variety} we mean a Zariski open and dense subset, $X$, of a (reduced) compact complex analytic space, $\overline X$, given together with the embedding of $X$ in $\overline X$.
By the Zariski topology on such an~$X$, we mean the one induced from the compactification~$\overline X$.
A holomorphic (respectively meromorphic) map between meromorphic varieties, $f:X\to Y$, will be called {\em definable} if it extends to a meromorphic map $\overline f:\overline X\to\overline Y$.
This is analogous to (and a generalisation of) the category of quasi-projective algebraic complex varieties with regular/rational morphisms.
As in the algebraic case, one could work more generally with ``abstract" meromorphic varieties -- so modelled on meromorphic varieties by finitely many charts and where the transitions functions are definable holomorphic maps -- but we will not do so here.

Let $\mathcal M$ be the structure where there is a sort for each irreducible meromorphic variety and a predicate for each Zariski closed subset of each finite cartesian product of sorts.
We denote the language of $\mathcal M$ by $L$, and by $\ccm$ the $L$-theory of $\mathcal M$.
It is a totally transcendental theory that admits elimination of quantifiers and imaginaries.
Note that all definable sets in~$\M$ are $0$-definable as every element of~$\M$ is named in the language~$L$, that is, $\dcl(\emptyset)=\M$.

Let $\AA$ denote the sort of the affine line.
Then $\AA(\M)=\CC$, the field structure on~$\CC$ is definable in~$\M$, and the full induced structure on~$\AA$ from~$\M$ is interdefinable (over parameters) with this field structure.
In this way, $\acf_0$ lives as a pure stably embedded reduct of $\ccm$.

Note that, as the terminology suggests, definable holomorphic/meromorphic maps between meromorphic varieties are definable in~$\ccm$.
While in algebraic geometry every regular map on a quasi-projective variety extends to a rational map on the projective closure, here the compactifiability condition (that the map extend to a meromorphic map on the given compactificiations) is not vacuous.
For example, complex exponentiation on the affine line does not extend to a meromorphic map on the projective line, and is indeed not definable in~$\ccm$.
In fact, a consequence of Chow's algebraicity theorem is that every definable holomorphic (respectively definable meromorphic) map on a quasi-projective complex algebraic variety is regular (respectively rational).

We work in a sufficiently saturated elementary extension $\overline\M\succeq\M$ and denote by $\C:=\AA(\overline\M)$ the corresponding (saturated) algebraically closed field extending~$\CC$.

Given an irreducible meromorphic variety, $X$, by a {\em generic} point of~$Y$ we mean an element $a\in X(\overline\M)$ not contained in $Y(\overline\M)$ for any proper Zariski closed subsets $Y\subset X$.
By quantifier elimination, this determines a complete type (over the empty set) in $\ccm$, the {\em generic type of~$X$}.
All complete types over the empty set arise this way: for any sort~$S$ of~$\overline\M$, and any point $a\in S(\overline\M)$, we have that~$a$ is a generic point of the {\em Zariski locus of~$a$}, denoted by $\loc(a)$, which is the smallest Zariski closed subset $X\subseteq S$ such that $a\in X(\overline\M)$.

\medskip
\subsection{The nonstandard Zariski topology}
Definable sets in~$\overline\M$ arise as {\em generic fibres} of families of definable sets in~$\M$; they are of the form $Z_b$ where~$Z$ is a definable subset of a product of irreducible meromorphic varieties $X\times Y$ projecting dominantly onto~$Y$, and $b\in Y(\overline\M)$ is generic.
We mean fibres, here, with respect to the second co-ordinate projection.
A definable subset~$D\subseteq X(\overline\M)$ is {\em Zariski closed} if it is of the form $Z_b$ as above with $Z\subseteq X\times Y$ Zariski closed.
To specify the parameters we say then that $D$ is {\em Zariski closed over~$b$}.
The Zariski topology on $X(\overline\M)$, thus defined, is noetherian. 
It agrees with the one induced from the Zariski topology on $\overline X(\overline\M)$.
See~\cite[$\S$2]{ret} for more details in the case when $X=\overline X$.

The following lemma -- which is a  consequence of the existence of a universal family of complex analytic subsets (the {\em Douady space}), and is implicit in the literature -- ensures that being $b$-definable and Zariski closed is equivalent to being Zariski closed over~$b$.

\begin{lemma}
\label{code}
Suppose $X$ is an irreducible meromorphic variety and $b$ is finite tuple from~$\overline\M$.
A $b$-definable Zariski closed subset of $X(\overline\M)$ is Zariski closed over~$b$.
\end{lemma}

\begin{proof}
Since the Zariski topology on $X(\overline\M)$ is induced by that on $\overline X(\overline\M)$, it suffices to show this in the case that $X=\overline X$.
Suppose $D\subseteq X(\overline\M)$ is Zariski closed.
This means that $D=Z_c$ where $Z\subseteq X\times Y$ is a Zariski closed subset, $Y$ is an irreducible compact complex analytic space, and $c\in Y(\overline\M)$ is a generic point.
We first claim that we can arrange things so that the fibres of $Z\to Y$ are pairwise distinct over a nonempty Zariski open subset of~$Y$.
Indeed, the theory of Douady spaces (along with Hironaka's flattening theorem) gives us a dominant meromorphic map $f:Y\to Y'$, where $Y'$ is an irreducible compact complex analytic space (living in the Douady space of $X$), such that for $y_1,y_2$ ranging in some nonempty Zariski open subset of~$Y$, we have $Z_{y_1}=Z_{y_2}$ if and only if $f(y_1)=f(y_2)$.
Moreover, there is a Zariski closed subset $Z'\subseteq X\times Y'$ (which is the restriction of the universal family of~$X$ to $Y'$) such that $Z_y=Z'_{f(y)}$ for all $y$ in some nonempty Zariski open subset of~$Y$.
See, for example, the second paragraph of the proof of Lemma~3.4 of~\cite{saturated}, for a detailed argument.
It follows that $D=Z_c=Z'_{f(c)}$, and $Z'\to Y'$ has the desired property that its fibres are pairwise distinct over a nonempty Zariski open subset of~$Z$.
So we may as well assume that $Z=Z'$ and $Y=Y'$.

Now assume that $D$ is $b$-definable.
We claim that $c\in\dcl(b)$.
Indeed, letting $\sigma\in\aut_b(\overline\M)$ be aribtrary, we have that
$$
Z_c
=
D
=
\sigma(D)
=
Z_{\sigma(c)}.$$
As the fibres of $Z\to Y$ are pairwise distinct over a nonempty Zariski open subset of~$Y$,
and $c,f(c)$ are both generic in $Y$, it follows that $\sigma(c)=c$.
So $c\in\dcl(b)$.

Finally, note that if~$D$ is Zariski closed over~$c$, and $c\in\dcl(b)$, then~$D$ is also Zariski closed over~$b$.
Indeed, letting $Y'=\loc(b)$, we get from quantifier elimination a dominant meromorphic map $f:Y'\to Y$ taking $b$ to $c$, so that if we take the base extension $Z':=Z_{Y'}\to Y'$ of $Z\to Y$ to $Y'$, then $D=Z'_b$ witnesses that $D$ is Zariski closed over $b$.
Here $Z_{Y'}$ is the fibred product in the meromorphic category: so the Zariski closure in $X\times Y'$ of the set of pairs $(x,y')$ such that~$y'$ is in the domain of~$f$ and $(x,f(y'))\in Z$.
\end{proof}

In the nonstandard setting ``irreducibility" depends on parameters.
We call a  nonstandard Zariski closed set~$D$ over~$b$, {\em irreducible over~$b$} if it cannot be written as the union of two proper nonstandard Zariski closed subsets over~$b$.
Equivalently, $D$ can be presented as $Z_b$ where $Z\subseteq X\times Y$ is Zariski closed and irreducible (for some $Y$ and $b\in Y(\overline\M)$ generic).
In this case, we have a unique {\em generic type of~$D$ over~$b$}, namely the type $p(x)\in S(b)$ saying that~$x$ is an element of~$D$ but not contained in any proper $b$-irreducible Zariski closed subset of~$D$ over~$b$.
(We say that $D$ is {\em absolutely irreducible} if it cannot be written as the union of two proper nonstandard Zariski closed subsets over any parameters, and this corresponds to asking that $Z\to Y$ being a {\em fibre space}, namely having irreducible fibres over some nonempty Zariski open subset of~$Y$.)

\medskip
\subsection{Internality and orthogonality in $\ccm$}
\label{subsect:intorth}
Let us recall the geometric characterisations of the stability theoretic notions around $\C$-internality in $\ccm$.
These results appear mostly in~\cite{ret}, where the reader can find further details.

First of all, for the purposes of geometric characterisations, we may restrict ourselves to types over finitely many parameters.
Indeed, by total transcendentality, the canonical base of a type is in the definable closure of a finite set.
Enumerating the finite set of parameters, and using that a finite cartesian product of sorts arises as a sort in its own right, we may further restrict ourselves to types of the form $p=\tp(a/b)$, where $a$ and $b$ are elements living in (possibly different sorts of)~$\overline\M$.
Letting $X:=\loc(a,b)$, $Y:=\loc(b)$, and $f:X\to Y$ the co-ordinate projection, we can realise $p$ as the generic type of the generic fibre over $f:X\to Y$.
In this way, complete types correspond to the bimeromorphic geometry of such maps.

\begin{fact}
\label{fact:gst-bg}
Fix a surjective meromorphic map $f:X\to Y$ between irreducible compact complex spaces, and let $p$ be the generic type of the generic fibre.
The following translates between the relevant model-theoretic and geometric notions:
\begin{itemize}
\item[(a)]
Stationarity of $p$ corresponds to $f$ being a {\em fibre space}, namely that the fibres of $f$ are irreducible over some nonempty Zariski open subset of~$Y$.
\item[(b)]
Weak $\C$-orthogonality of $p$ corresponds to the condition that the induced map $f^*:\CC(Y)\to \CC(X)$ on meromorphic function fields is an isomorphism.
\item[(c)]
Suppose $p$ is stationary.
Internality of~$p$ to~$\C$ corresponds to some base extension\footnote{To be more accurate, one replaces the fibre product $X\times_YZ$ with its unique irreducible component projecting onto~$Z$. That there is a unique such irreducible component uses that~$f$ is a fibre space.}
of~$f$,
$$f_Z:X\times_YZ\to Z,$$
admitting a meromophic embedding into $Z\times \PP^n$ over~$Z$, for some $n>0$.
Here $Z\to Y$ is an irreducible compact complex space with a surjective meromorphic map to $Y$.
\item[(d)]
Suppose $p$ is stationary.
Fundamental $\C$-internality of~$p$ corresponds to when $Z\to Y$ in part~(c) is just $f:X\to Y$ itself.
\end{itemize}
\end{fact}

\begin{proof}
Part~(a) is~\cite[Lemmas~2.7 and~2.11]{ret}.

Part~(b), is well known --  see, for example, \cite[Remark~7.7]{abred} -- but we sketch the idea here.
First of all, given a generic point $a\in X(\overline\M)$, there is a natural identification of $\CC(X)$ with $\dcl(a)\cap\C$,  the identification is given by evaluating at~$a$, namely $g\mapsto g(a)$.
If we fix $a\models p$, then~$a$ is generic in~$X$ and $f(a)$ is generic in ~$Y$, so that $f^*$ corresponds to the containment $\dcl(f(a))\cap\C\subseteq\dcl(a)\cap\C$.
So, setting $E:=\dcl(a)\cap\C$, we need to show that $\tp(a/f(a))$ isolates $\tp(a/f(a)\C)$ if and only if $E\subseteq\dcl(f(a))$.
The left-to-right direction is clear, and for the right-to-left direction we use that, as~$\C$ is definably closed and stable embedded, we do have that $\tp(a/f(a)E)$ isolates $\tp(a/f(a)\C)$.

Part~(c) is~\cite[Proposition~4.4]{ret}.
The way $Z$ is obtained under the assumption that $p$ is $\C$-internal is as the locus of the additional parameters used to witness internality: that is, if $a\models p$ and $b$ are such that $a\ind_{f(a)} b$ and $a\in\dcl(f(a)b\C)$, then we can take $Z=\loc(f(a),b)$.

Part~(d) follows from the above description of~$Z$ witnessing internality, noting that $p$ is fundamental $\C$-internal if and only if there are $a,b$ realising $p$, with $a\ind_{d} b$ and $a\in\dcl(db\C)$ where $d=f(a)=f(b)$.
\end{proof}

It is maybe worth specialising to the case of types over the empty set; so this is the absolute case when $Y$ is zero-dimensional.
That is, $p$ is the generic type of an irreducible compact complex space $X$.
In that case, part~(a) of Fact~\ref{fact:gst-bg} says that $p$ is always stationary (this is because $\dcl(\emptyset)=\M$ is a model), part~(b) says that $p$ is weakly $\C$-orthogonal if and only if $X$ has no nonconstant meromorphic functions, 
and part~(c) says that $p$ is $\C$-internal if and only if $X$ is bimeromorphic to a projective algebraic variety, i.e., $X$ is {\em Moishezon}.
Note that $X$ is Moishezon if and only if $\trdeg_{\CC}(\CC(X))=\dim X$.

\medskip
\subsection{The K\"ahler case}
Sometimes we assume that the the compact complex spaces~$X$ and~$Y$, where~$p$ is the generic type of the generic fibre of $f:X\to Y$, are of K\"ahler-type in the sense of Fujiki~\cite{fujiki-kaehler}.
We will say, in that case, that {\em $p$ is of K\"ahler-type.}
Such an assumption is motivated, model-theoretically, by ``essential saturation" which often gives automatic uniformity; see~\cite{saturated} for details.
For example, under the K\"ahler-type assumption, $\C$-internality of~$p$ is equivalent to asking that over a nonempty Zariski open subset of $Y$ the fibres of $f$ are Moishezon (compare with condition~(c) of Fact~\ref{fact:gst-bg}).

\medskip
\subsection{Algebraic dimension and reduction}
Given an irreducible meromorphic variety~$X$, by the {\em algebraic dimension} of~$X$, denoted by $a(X)$,
we mean the transcendence degree of $\CC(\overline X)$, the field of meromorphic functions on the compactification.
In particular, $a(X)=\dim(X)$ if and only if~$\overline X$ is Moishezon, and $a(X)=0$ if and only if it admits no definable nonconstant meromorphic functions.
Stated model-theoretically, as a consequence of Fact~\ref{fact:gst-bg}, the former is equivalent to the generic type of~$X$ being $\C$-internal, and the latter is equivalent to the generic type of~$X$ being $\C$-orthogonal.
(We use here that $\dcl(\emptyset)=\M$ is a model.)

An {\em algebraic reduction} of an irreducible compact complex space~$X$ is a fibre space $f:X\to Y$ such that $Y$ is Moishezon and every surjective meromorphic map from $X$ to a Moishezon space factors through $f$.
This forces $f^*:\CC(Y)\to\CC(X)$ to be an isomorphism, so that $a(X)=\dim(Y)$.
Algebraic reductions always exist and are unique up to bimeromorphism, see~\cite[Chapter~VII]{scv7} for more details.
By Fact~\ref{fact:gst-bg}, the generic type of the generic fibre of~$f$ is weakly $\C$-orthogonal.
In fact, this is how all weakly $\C$-orthogonal types over parameters from~$\C$ arise.

\medskip
\subsection{Linear spaces}
\label{subsect:linear}
One source of types internal to~$\C$ are linear spaces.
Suppose $Y$ is an irreducible compact complex space~$Y$.
A {\em linear space over $Y$}, $f:L\to Y$, is a $(Y\times\mathbb C)$-module in the category of complex spaces over $Y$.
That is, there are holomorphic maps
$$
\xymatrix{
L\times_YL\ar[dr]_{f\pi_2}\ar[rr]^{+}&&L\ar[dl]^{f}\\
&Y&
}
$$
and 
$$
\xymatrix{
(Y\times\mathbb C)\times_YL\ar[dr]_{f\pi_2}\ar[rr]^{\lambda}&&L\ar[dl]^{f}\\
&Y&
}
$$
and a {\em zero section} $z:Y\to L$ to $f$,
satisfying the usual module axioms.
In particular, each fibre of $f$ is endowed with the structure of a finite dimensional $\mathbb C$-vector space, uniformly over~$Y$.
The dimension of this vector space is $\dim L-\dim Y$, which we call the {\em generic rank} of the linear space.
See~\cite[$\S1.4$]{fischer76} for details.

Another perspective on linear spaces is that of $\mathcal O_Y$-modules.
Indeed, taking the sheaf of linear forms on a linear space over~$Y$, namely the sheaf of homomorphisms to $Y\times\mathbb C$ over~$Y$, gives us a coherent  $\mathcal O_Y$-module.
That this defines an antiequivalence of categories is explained in~\cite[$\S1.6$]{fischer76}.

There is a natural compactification, $\overline L$, of $L$ making it an irreducible meromorphic variety (and hence a sort of $\mathcal M$), and such that the linear structure (namely $f, +,\lambda, z$) extends meromorphically to $\overline L\to Y$, making it definable in $\mathcal M$.
Indeed, $\overline L$ can be taken to be the {\em projective linear space over~$Y$} corresponding to $L\times\CC\to Y$; so
this is just the relativisation of the usual compactification of $\mathbb C^n$ as $\mathbb P_n(\mathbb C)$.
See~\cite[$\S1.9$]{fischer76} for details.

Linear spaces give rise to $\C$-internality.
Indeed, interpreting $L\to Y$ in $\overline{\mathcal M}$, the fibres are uniformly definable vector spaces over~$\mathcal C$.
In particular, for $b\in Y(\overline\M)$ generic, $L_b$ admits a $b$-definable $\mathcal C$-vector space structure, and hence is definably isomorphic to $\C^n$ over any choice of basis.

It seems plausible that all definable vector spaces over~$\C$ arise in this way:

\begin{question}
\label{question}
Is every $b$-definable $\C$-vector space in~$\overline\M$ $b$-definably isomorphic to~$L_b$ for some linear space $L\to \loc(b)$?
\end{question}

\medskip
\subsection{Moishezon maps}
More generally, a surjective meromorphic map $f:X\to Y$ is called {\em Moishezon} if it factors meromorphically over $Y$ into a projective linear space over~$Y$.
These too, then, give rise to $\C$-internal types.
For a more detailed treatment of the model-theoretic interpretation of linear spaces and Moishezon maps, see~\cite[$\S$3]{ret}.

When $X$ is not Moishezon, the algebraic reduction cannot be a Moishezon map: if $Y$ is Moishezon and $f:X\to Y$ is Moishezon, then $X$ itself is Moishezon. 
However, it can still happen that the generic type of a generic fibre of an algebraic reduction $f:X\to Y$ is $\C$-internal.
Model theoretically, this corresponds to the generic type of $X$ (over the emptyset) being analysable in~$\C$ in two steps.
So this gives rise to examples of $\C$-internality (with weak $\C$-orthogonality) that do {\em not} arise from Moishezon maps.

\bigskip
\section{An application to principal $G$-bundles}
\label{sec:bundle}

\noindent
The following is the natural relativisation of a uniquely transitive action of an algebraic group on a meromorphic variety:

\begin{definition}
Suppose $G$ is an algebraic group over~$\CC$.
By a {\em principal meromorphic $G$-bundle} we mean a definable surjective holomorphic map $f:X\to Y$, where $X$ and $Y$ are irreducible meromorphic varieties, together with a definable holomorphic map $X\times G(\CC)\to X$ over~$Y$, which restricts to a uniquely transitive right action of $G(\CC)$ on each fibre of~$f$.
\end{definition}

\begin{remark}
\begin{itemize}
\item[(a)]
By the Fischer-Grauert  theorem~\cite{fischergrauert}, after possibly shrinking~$Y$, $f:X\to Y$ will be a holomorphic fibre bundle (in the sense that it is locally trivial in the holomorphic category) with structure group~$G$.
Our assumption that~$f$ and the action of~$G$ extend meromorphically to compactifications (this is the content of ``definable holomorphic") means that the fibre bundle has {\em meromorphic structure} in the sense of Fujiki~\cite{fujiki-fiber}.
\item[(b)]
While we will not do so, it does make sense to talk about principal meromorphic bundles for arbitrary meromorphic groups (in the sense of~\cite{mergroup}), not just for algebraic groups.
\end{itemize}
\end{remark}

Passing to the nonstandard model~$\overline\M$, and taking $b\in Y(\overline\M)$ generic, we have that the generic fibre, $X_b$, of a principal meromorphic $G$-bundle, is a definable (right) torsor for $G(\C)$.
In particular, $X_b$ is $\C$-internal.
Our goal in this section is to investigate the relationship between $G$ and the binding group for $X_b$ relative to~$\C$.
This will allow us to apply the general model theory of binding groups to deduce structural results about principal meromorphic bundles (namely, Corollary~\ref{cor:tt}).

We will be focusing on principal meromorphic bundles where the base is of algebraic dimension zero.

\begin{lemma}
\label{gbundle-bg}
Suppose
$G$ is an algebraic group over~$\CC$, and $X\to Y$ is a principal meromorphic $G$-bundle with $a(Y)=0$.
Let $b\in Y(\overline\M)$ be generic, $M\preceq\overline\M$ a prime model over~$b$,
$a\in X_b(M)$, and $q:=\tp(a/b)$.
Then
\begin{itemize}
\item[(a)]
$q$ is fundamental ${\mathcal C}$-internal and weakly orthogonal to ${\mathcal C}$,
\item[(b)]
the restriction map $\aut_b(X_b/\C)\to \aut_b(q/\C)$ is an isomorphism, and
\item[(c)]
the dual binding group of~$q$ relative to~$\C$ over~$b$ is of the form $H(\C)$ where $H\leq G$ is an algebraic subgroup defined over~$\CC$, with right action on $q(\overline\M)$ induced by the given right action of~$G$ on $X_{b}$. 
\end{itemize}
In particular, the binding group of the generic fibre of $X\to Y$ is definably isomorphic to the $\C$-points of a complex algebraic subgroup of~$G$.
\end{lemma}

\begin{proof}
(a).
As $M$ is prime over~$b$, we have that 
$M\cap\C$ is algebraic over $\dcl(b)\cap\C$.
But the assumption that $a(Y)=0$ implies that $\dcl(b)\cap\C=\CC$.
In particular, $\dcl(ab)\cap\C=\CC$, so that, by Lemma~\ref{lem:wochar}, $q=\tp(a/b)$ is weakly $\C$-orthogonal.
That it is fundamental $\C$-internal follows from the fact that $G(\C)$ acts transitively on $X_b$.

(b).
Restriction is a surjective homomorphism.
It is injective as $X_b\subseteq\dcl(ab\C)$.
 
(c).
For any $a'\models q$ there is unique $g\in G({\mathcal C})$ with $a' = ag$.  The construction in  the proof of Theorem~\ref{bgt} therefore produces the dual binding group as a $b$-definable subgroup of $G({\mathcal C})$.
Any such is of the form $H(\C)$ for some algebraic subgroup $H\leq G$ defined over  $\dcl(b)\cap {\mathcal C} = {\mathbb C}$.

The ``in particular" clause follows from putting~(ii) and~(iii) together with the fact that the binding and dual binding groups are definably isomorphic.
\end{proof}

\begin{proposition}
\label{fullbg}
Suppose $X\to Y$ is a principal meromorphic $G$-bundle with $a(Y)=0$.
Let $b\in Y(\overline\M)$ be generic.
Then the following are equivalent:
\begin{itemize}
\item[(i)]
$X_b$ is a complete type over~$b$.
\item[(ii)]
If $H$ is a proper algebraic subgroup of~$G$ over~$\CC$, and $a\in X_b(\overline\M)$, then $aH(\C)$ is not defined over~$b$.
\item[(iii)]
For $M\preceq\overline\M$ a prime model over~$b$, and
$a\in X_b(M)$, the dual binding group of $q:=\tp(a/b)$, under the identification of Lemma~\ref{gbundle-bg}(c), is all of $G(\C)$.
\item[(iv)]
The generic type of $X_b$ over~$b$ is weakly $\C$-orthogonal.
\end{itemize}
\end{proposition}

\begin{proof}
(i)$\implies$(ii).
Suppose $X_b$ is a complete type over~$b$.
Fix $H\leq G$ an algebraic subgroup over~$\CC$, and $a\in X_b(\overline\M)$.
If $aH(\C)$ is defined over~$b$ then membership in it is part of $\tp(a/b)$, which is realised by every element of $X_b$.
Hence $aH(\C)=X_b$.
But the unique-transitivity of the action of $G(\C)$ on $X_b$ would then force $H=G$.

(ii)$\implies$(iii).
Lemma~\ref{gbundle-bg} gives us that~$q$ is $\C$-internal and weakly $\C$-orthogonal with dual binding group of the form $H(\C)$, for some algebraic subgroup $H\leq G$, with the induced right action.
So $aH(\C)=q(\overline\M)$, from which it follows that $a H(\C)$ is $b$-definable.
By~(ii), $H=G$, as desired.

(iii)$\implies$(iv).
Let~$q$ be as in condition~(iii).
As $G(\C)$ acts transitively on~$X_b$, and the dual binding group acts transitively on $q(\overline\M)$, condition~(iii) tells us that $q(\overline\M)=X_b(\overline\M)$.
In particular, $q$ is the generic type of $X_b$, and it is weakly $\C$-orthogonal by Lemma~\ref{gbundle-bg}(a).

(iv)$\implies$(i).
Since $X$ is irreducible, $X_b$ is irreducible over~$b$, and hence has a unique generic type, say $r=\tp(e/b)$.
Let $N\preceq\overline\M$ be a prime model over~$e$, so that $N\cap\C$ is algebraic over $\dcl(e)\cap\C$.
That~$r$ is weakly orthogonal to~$\C$ implies that
$\dcl(e)\cap\dcl(\C b)=\dcl(b)$, by Lemma~\ref{lem:wochar}.
Hence $\dcl(e)\cap \C=\dcl(b)\cap \C=\CC$ as $a(Y)=0$.
That is, $N\cap\C=\CC$.
Now, suppose $a\in X_b(N)$.
As $G(N)=G(\CC)$ acts transitively on $X_b(N)$, there is $g\in G(\CC)$ such that $ga=e$.
But the action by~$g$, since~$g$ is $0$-definable, is a definable bijection between $\tp(a/b)$ and $r=\tp(e/b)$.
It follows that $\tp(a/b)$ is also generic in $X_b$, and hence $tp(a/b)=r$.
That is, every $N$-point of $X_b$ realises~$r$.
Since $N\preceq\overline\M$, it follows that $X_b(\overline\M)=r(\overline\M)$.
\end{proof}

We deduce the following characterisation for what might be called ``strong irreducibility" of a principal meromorphic bundle -- following the language for vector bundles introduced in~\cite[$\S$1]{toma}.
But first let us make precise what we mean by sub-bundles of principal meromorphic $G$-bundles.

\begin{definition}
Suppose $G$ is an algebraic group over~$\CC$ and $X\to Y$ is a principal meromorphic $G$-bundle.
By a {\em principal meromorphic sub-bundle} of $X\to Y$ we mean a Zariski closed subset $Z\subseteq X$ and a nonempty Zariski open subset $U\subseteq Y$ such that for some algebraic subgroup $H\leq G$ over~$\CC$, the action of $G(\CC)$ on~$X$ induces an action of $H(\CC)$ on~$Z|_U$ making $Z|_U\to U$ a principal meromorphic $H$-bundle.
\end{definition}

\begin{corollary}
\label{cor:tt}
Suppose $G$ is a complex algebraic group, $f:X\to Y$ is a principal meromorphic $G$-bundle, and $a(Y)=0$.
Then the following are equivalent:
\begin{itemize}
\item[(i)]
$X$ admits no proper Zariski closed subsets projecting dominantly onto~$Y$.
\item[(ii)]
$f:X\to Y$ admits no proper principal meromorphic sub-bundles.
\end{itemize}
\end{corollary}

\begin{proof}
Fix $b\in Y(\overline\M)$ generic.

Note first of all that, for $Z\subseteq X$ Zariski closed and projecting dominantly onto~$Y$, we have that $Z=X$ if and only if $Z_b=X_b$.
This is by irreducibility of~$X$.
Hence, by quantifier elimination for~$\ccm$, condition~\ref{cor:tt}(i) is equivalent to condition~\ref{fullbg}(i).
We will show, and this will suffice, that condition~\ref{cor:tt}(ii) is equivalent to condition~\ref{fullbg}(ii).

That~\ref{fullbg}(ii) implies~\ref{cor:tt}(ii) is clear:
If $Z\to X$ is a principal meromorphic sub-bundle associated to $H\leq G$, and $a\in Z_b(\overline\M)$, then $aH(\C)=Z_b$ is $b$-definable, so that $H=G$ by~\ref{fullbg}(ii), and hence $Z_b=X_b$, implying $Z=X$.

For the converse, we assume~\ref{cor:tt}(ii) and verify~\ref{fullbg}(ii).
Suppose $H\leq G$ is an algebraic subgroup over~$\CC$, and $a\in X_b(\overline\M)$ has a $b$-definable orbit $aH(\C)$.
We first argue that  $aH(\C)$ is Zariski closed over~$a$ in $X(\overline\M)$.
To do so we need to exhibit $aH(\C)$ as the fibre over~$a$ of some Zariski closed subset of $X\times X$.
Indeed, consider the fibred product $\widetilde X=X\times_YX$, as a Zariski closed subset of $X\times X$, and equipped with the first co-ordinate projection $\widetilde X\to X$.
Note that there is a definable biholomorphism:
$$\xymatrix{
X\times G(\CC)\ar[dr] \ar[rr]^{\rho}_{\approx}&& \widetilde X\ar[dl]\\
&X&
}$$
given by $\rho(x,g):=(x,xg)$.
As $H(\CC)$ is Zariski closed in $G(\CC)$, we have that $\rho(X\times H(\CC))$ is Zariski closed in~$\widetilde X$, and hence in $X\times X$.
As $\rho(X\times H(\CC))_a=aH(\C)$, this witnesses that $aH(\C)$ is Zariski closed over~$a$.

So the orbit $aH(\C)$ is both $b$-definable and Zariski closed.
Lemma~\ref{code} now ensures that $aH(\C)$ is Zariski closed over~$b$.
That is, $aH(\C)=Z_b$ for some Zariski closed $Z\subseteq X$.
As the induced action of $H(\C)$ is uniquely transitive on $Z_b$, we have, back in the standard model, that for some nonempty Zariski open $U\subseteq Y$, there is an induced action of $H(\CC)$ on $Z|_U$ over~$U$ such that $Z|_U\to U$ is a principal meromorphic $H$-bundle.
That is, $Z$ is a principal meromorphic sub-bundle of~$f$.
Hence, by~\ref{cor:tt}(ii), we must have $Z=X$, and hence $aH(\C)=X_b$, which forces $H=G$, as desired.
\end{proof}

\bigskip
\section{$\ccm$-Galois groups}
\label{sec:ccmgal}

\noindent
We continue to work in a sufficiently saturated $\overline\M\models\ccm$ with $\C=\AA(\overline\M)$ the corresponding (saturated) algebraically closed field extending~$\CC$.

The following points out that binding groups relative to~$\C$ are algebraic.

\begin{lemma}
\label{bgag}
Suppose $\Phi$ is a $\C$-internal (partial) type with binding group~$\G$ relative to~$\C$.
Then $\G$ is definably isomorphic, over possibly additional parameters, to $G(\C)$ for some algebraic group $G$ over~$\C$.
\end{lemma}

\begin{proof}
We may assume that~$\Phi=p\in S(B)$ is complete, fundamental $\C$-internal, and weakly $\C$-orthogonal.
Indeed, Proposition~\ref{prop:redtofundwo} tells us that a binding group for~$\Phi$  is definably isomorphic to the binding group of a complete type with those properties.
(This already requires working over a potentially larger set of parameters.)
But now apply Theorem~\ref{bgt} and we have a dual binding group $\widehat\G$ that is definably isomorphic to $\G$ (over a realisation of~$p$), and such that
$\widehat\G\subseteq\dcl(B\C)$.
As~$\C$ is stably embedded in $\overline\M$ and the full induced structure is that of a pure field, and as algebraically closed fields admit elimination of imaginaries, we get that $\widehat\G$ is $B$-definably isomorphic to a definable group in $(\C,+,\times)$.
It follows that $\widehat\G$, and hence~$\G$, is definably isomorphic to $G(\C)$ for some algebraic group~$G$.
\end{proof}

Inspecting the proof, we obtain the following control on the parameters if we happen to start with a complete, fundamental $\C$-internal, and weakly $\C$-orthogonal type, and we focus instead on the dual binding group:

\begin{corollary}
\label{cor:bgag}
If $p\in S(B)$ is fundamental $\C$-internal and weakly $\C$-orthogonal, and $\widehat\G$ is the dual binding group for~$p$ relative to~$\C$, then $\widehat\G$ is $B$-definably isomorphic to $G(\C)$ where~$G$ is an algebraic group over~$\dcl(B)\cap\C$.
\end{corollary}

\begin{proof}
This was explicitly established in the proof of Lemma~\ref{bgag}.
\end{proof}

\begin{definition}
\label{def:galois}
We say that an algebraic group~$G$ over $\C$ is a {\em $\ccm$-Galois group} if 
$G(\C)$ is definably isomorphic to the binding group of some (partial) $\C$-internal type.
More precisely, there is a parameter set~$B$ and a (partial) $\C$-internal type~$\Phi$ over~$B$, with a definable binding group~$\G$ relative to~$\C$ over~$B$, such that $G(\C)$ and~$\G$ are definably isomorphic over possibly additional parameters.
\end{definition}

Being a $\ccm$-Galois group is witnessed more explicitly as follows:

\begin{lemma}
\label{galois-dual}
Suppose $G$ is a $\ccm$-Galois group.
Then there exists:
\begin{itemize}
\item
$p\in S(B)$ fundamental $\C$-internal and weakly $\C$-orthogonal, and
\item
$H$ an algebraic group over $\dcl(B)\cap\C$,
\end{itemize}
such that $G$ and $H$ are isomorphic algebraic groups and $H(\C)$ is the dual binding group for~$p$ relative to~$\C$.
\end{lemma}

\begin{proof}
Because of Proposition~\ref{prop:redtofundwo}, we can always take~$\Phi$ in Definition~\ref{def:galois}, witnessing that $G$ is a $\ccm$-Galois group, to be a fundamental $\C$-internal and weakly $\C$-orthogonal type $p\in S(B)$, for some choice of~$B$.
And, because of Corollary~\ref{cor:bgag}, the dual binding group of such~$p$ is of the form $H(\C)$ for some algebraic group $H$ over $\dcl(B)\cap\C$.
Since the binding group is definably isomorphic to the dual binding group, $G$ and $H$ are isomorphic as algebraic groups (over~$\C$).
\end{proof}

In Definition~\ref{def:galois} we allowed the parameters to vary.
It may also be of interest to fix parameters~$B$ and talk about $\ccm$-Galois groups {\em over $B$}, meaning that the conclusions of Lemma~\ref{galois-dual} hold with this~$B$.

The question of which algebraic groups are $\ccm$-Galois groups was raised as Problem~4.12 in~\cite{jjp}.
It was shown there, based on a construction of Toma~\cite{toma}, that $\PGL_2$ is a $\ccm$-Galois group.
On the other hand, in~\cite[Proposition~7.1]{abred}, it was shown that the binding group of a type over parameters in~$\C$ cannot be linear algebraic.
Little else seems to have been said on the subject, and our goal here is to collect some further examples and observations.

\medskip
\subsection{Simple abelian varieties as $\ccm$-Galois groups}
\label{subsect:abelian}
We show that every simple abelian variety is a $\ccm$-Galois group.
In Appendix~\ref{campana-appendix}, Fr\'ed\'eric Campana shows that if we start with a simple abelian variety~$A$, and any other abelian variety~$B$, then there is a complex torus $T$ that is an extension of~$B$ by~$A$, and such that the quotient map $f:T\to B$ is the algebraic reduction of~$T$.
In particular, the following two properties hold of the fibration $f:T\to B$.
\begin{itemize}
\item[(i)]
$f^*:\CC(B)\to\CC(T)$ is an isomorphism of the meromorphic function fields,
\item[(ii)]
for generic $b\in B(\overline\M)$, the fibre $T_b$ is definably isomorphic (over additional parameters) to $A(\C)$.
\end{itemize}
Fix $b\in B(\overline\M)$ generic, and let $p\in S(b)$ be the generic type of the fibre $T_b$.
Condition~(i) says precisely that $p$ is weakly orthogonal to~$\C$, see Fact~\ref{fact:gst-bg}(b).
Condition~(ii) implies, in particular, that $p$ is $\C$-internal.
Let $\G$ be the binding group of~$p$.
By Lemma~\ref{bgag}, $\G$ is definably isomorphic to $G(\C)$ for some algebraic group~$G$ over~$\C$.
Isolation and generecity of~$p$ imply that $p(\overline\M)$ contains a dense Zariski open subset of $T_b$, and the latter is definably isomorphic to $A(\C)$.
The definable action of~$\G$ on $p(\overline\M)$ thus induces an algebraic action of~$G$ on~$A$ with a Zariski dense orbit.
This forces~$G$ to be isomorphic to~$A$ (with the left multiplication action of~$A$ on itself).
Hence, $A$ is a $\ccm$-Galois group, as claimed.

\medskip
\subsection{The multiplicative group as a $\ccm$-Galois group}\label{subsect:gm}
Suppose we have a linear space $f:L\to Y$ of generic rank~$1$, and such that  $f^*:\CC(Y)\to\CC(\overline L)$ is an isomorphism of the meromorphic function fields.
(See Section~\ref{subsect:linear} on linear spaces.)
Such linear spaces do exist: any nontrivial line bundle on a two-dimensional complex torus without curves would be an example.
We will show that~$f$ witnesses that the multiplicative group $\mathbb G_m$ is a $\ccm$-Galois group.

Fix $b\in Y(\overline\M)$ generic, and let $p\in S(b)$ be the generic type of the fibre $L_b$.
Again, the condition that $f^*$ is an isomorphism of meromorphic function fields ensures that~$p$ is weakly $\C$-orthogonal.
On the other hand, $L_b$ having a $B$-definable $\mathcal C$-vector space structure ensures that $p$ is $\C$-internal.
In fact, $p$ is fundamental $\mathcal C$-internal:
Suppose $a\models p$. 
Then $a\neq z(b)$, and hence $a$ is a basis for the ($1$-dimensional) $\C$-vector space $L_b$, and we obtain a $ba$-definable isomorphism $\mu:L_b\to \C$.

It follows from isolation and genericity of~$p$ that $\mu$ restricts to a definable bijection between $p(\overline\M)$ and a cofinite subset $Q\subseteq\C$.
Letting $G$ be an algebraic group over~$\C$ such that $G(\C)$ is definably isomorphic to the binding group of~$p$, we get a uniquely transitive algebraic action of $G$ on the affine line minus a finite set of points.
The only possibilities for such are the canonical actions of additive or multiplicative groups.
But as $Q\neq\C$ (since $p(\overline\M)\neq L_b$, since $z(b)\not\models p$), we can't be looking at the additive group acting on the affine line.
Hence $G=\mathbb G_m$, as desired.

\medskip
\subsection{A family of linear $\ccm$-Galois groups}
\label{sect:pgl}
We exhibit here, for each $n\geq 1$, $\ccm$-Galois groups of dimension at least $2n-1$ that are algebraic subgroups of $\PGL_{2n}$.
An earlier version of this paper claimed that $\PGL_{2n}$ itself is a $\ccm$-Galois group, but an error in our argument was pointed out to us by Remi Jaoui.

It is the case that $\PGL_2$ is a $\ccm$-Galois group.
In~\cite{toma} Toma shows that there are compact complex surfaces~$Y$ with $a(Y)=0$ that admit rank~$1$ projective linear spaces, $P\to Y$, such that no proper Zariski closed subset of~$P$ projects onto~$Y$.
A consequence of this, established in~\cite[$\S$4.3]{jjp}, is that, for generic $b\in Y(\overline\M)$,  the binding group~$\G$ of $P_b$ relative to~$\C$ over~$b$ is definably isomorphic to $\PGL_2(\C)$.
Our goal here is to extend this argument, as much as possible, to higher dimensions.
Toma's construction does extend: 
in Appendix~\ref{toma-appendix} below, he shows that, for all $n\geq 1$, there exist irreducible compact complex spaces~$Y$ of dimension $2n$ that admit projective linear spaces $P\to Y$ of rank $2n-1$, such that $a(P)=0$ and~$P$ contains no proper Zariski closed subsets that project onto~$Y$.
These examples arise by taking~$Y$ to be a compact irreducible hyperk\"ahler manifold with trivial Picard group, and~$P\to Y$ the projectivised tangent bundle of~$Y$.

Let $b\in Y(\overline\M)$ be generic, and let~$\G$ be the binding group of $P_b$ relative to~$\C$ over~$b$.
We know that $P_b$ is definably isomorphic to $\PP^{2n-1}(\C)$, over possibly additional parameters, and that this induces an action of~$\G$ on $\PP^{2n-1}(\C)$ by definable bijections.
In fact, as we now explain, there is an identification of $P_b$ with $\PP^{2n-1}(\C)$ such that the induced action is by automorphisms.
First of all, every projective linear space trivialises after some base extension.
That is, there is an irreducible meromorphic variety~$V$, with a dominant definable holomorphic map $V\to Y$, and a definable isomorphism $f:P\times_YV\to V\times\PP^{2n-1}$ over~$V$. (See, for example, \cite[Lemma~3.3]{ret}.)
Letting $v\in V(\overline\M)$ lie over~$b\in Y(\overline\M)$, we have that
$$f_v:P_b\to\PP^{2n-1}(\C)$$
is a (nonstandard) isomorphism.
We claim that this identification induces an action of~$\G$ on $\PP^{2n-1}(\C)$ by automorphisms.
Indeed, for any $g\in\G$, the induced action of~$g$ on $\PP^{2n-1}(\C)$ is given by $f_v\circ g\circ f_v^{-1}$.
But, if we let $\sigma\in\aut_{\C b}(\overline\M)$ act on $P_b$ as~$g$ does, then
$f_{\sigma(v)}:P_b\to\PP^{2n-1}(\C)$ is another isomorphism, and
$$f_v\circ g\circ f_v^{-1}=f_v\circ f_{\sigma(v)}^{-1}:\PP^{2n-1}(\C)\to\PP^{2n-1}(\C)$$
is an automorphism.
What we obtain is a definable embedding of~$\G$ in $\PGL_{2n}(\C)$.
The image is therefore of the form $G(\C)$ where~$G$ is an algebraic subgroup of $\PGL_{2n}$ that is therefore a $\ccm$-Galois group.

That~$P$ has no proper Zariski closed subsets projecting onto~$Y$ ensures that $P_b$ is the set of realisations of a complete type~$p$ over~$b$.
That $a(P)=0$ forces~$p$ to be weakly orthogonal to~$\C$.
Hence~$\G$ acts transitively on $P_b$, and so $\dim (G)\geq 2n-1$.

Computing precisely what the $\ccm$-Galois group is in these cases, namely when~$Y$ is a compact irreducible hyperk\"ahler manifold with trivial Picard group and~$P\to Y$ is the projectivised tangent bundle, would certainly be of interest.

\bigskip
\section{Characterising the linear case}
\label{sect:linalg}

\noindent
We have seen that $\mathbb G_m$ and certain algebraic subgroups of $\PGL_{2n}$ are $\ccm$-Galois groups.
In this section we investigate further the linear case.
In particular, we ask which $\C$-internal types have binding group linear algebraic?

Here is one answer that takes a little further observations already appearing in~\cite{abred}.
It is in terms of Fujiki's~\cite{fujiki-alb} relative Albanese in the K\"ahler-type setting.

\begin{proposition}
\label{prop:linalgbg}
Suppose $p$ is the generic type of the generic fibre of a fibre space
$f:X\to Y$ between K\"ahler-type compact complex spaces.
Suppose $p$ is $\C$-internal and weakly $\C$-orthogonal.
Then the binding group of~$p$ is definably isomorphic to the $\C$-points of a linear algebraic group if and only if the relative Albanese of $f$ is trivial.
Moreover, in this case $f$ is Moishezon.
\end{proposition}

\begin{proof}
As pointed out in the proof of~\cite[Proposition~7.1]{abred}, the assumptions that~$X$ is K\"ahler-type, that~$f$ is a fibration with general fibre\footnote{Here, by a {\em general fibre} having property~$P$ we mean that the fibres over some dense Zariski open subset of~$Y$ have property~$P$.} Moishezon (which follows from $\C$-internality of~$p$), allow us to invoke a theorem of Fujiki's~\cite[Theorem 2]{fujiki-alb} asserting that a relative Albanese exists.
Namely, there is a fibration $\operatorname{Alb}(f)\to Y$ whose general fibre is an abelian variety, and a meromorphic map $\alpha:X\to\operatorname{Alb}(f)$ over $Y$ which is the Albanese on the general fibre of~$f$.
Moreover, Fujiki shows that $\alpha$ is Moishezon.

If the binding group of~$p$ is linear then the general fibre of $f$ admits the action of a linear algebraic group with Zariski dense orbit.
But this implies that the Albanese of the general fibre is trivial, so that $\operatorname{Alb}(f)=Y$ and $\alpha=f$.
That is, the relative Albanese of $X$ over $Y$ is trivial and $f$ is Moishezon, as desired.

For the converse, suppose the binding group~$\G$ of~$p$ is not linear.
Let us work over $B=\acl(b)$ where $b\in Y(\overline{\mathcal M})$ is generic.
The main results of~\cite{abred} (see Theorems~3.7 and~4.2 there) give us an {\em abelian reduction} $p\to p_{\operatorname{ab}}$ over~$B$ where the binding group of $p_{\operatorname{ab}}$ is the abelian part of  $\G$, say~$\A$.
As $\G$ is not linear, $\A$ is not trivial.
Geometrically, $p\to p_{\operatorname{ab}}$ arises from a commutative diagram of surjective meromorphic maps
$$
\xymatrix{
X\ar[dr]_f\ar[rr]&&Z\ar[dl]^g\\
&Y
}$$
where 
$p_{\operatorname{ab}}$ is the generic type of the fibre over~$b$ of the fibration  $g:Z\to Y$.
Note that $\dim Z>\dim Y$ as $p_{\operatorname{ab}}$ is nonalgebraic (since its binding group, $\A$, is not trivial).
Since $\A$ is definably isomorphic to the $\C$-points of an abelian variety, the general fibre of $g$ admits the action of an abelian variety with Zariski dense orbit.
It follows that the general fibre of $g$ is itself a (positive dimensional) abelian variety.
In particular, the commuting diagram implies that the general fibre of $f$ has nontrivial Albanese, and hence $\operatorname{Alb}(f)$ is nontrivial.
\end{proof}

\begin{remark}
It follows that if $f:X\to Y$ is as in the above proposition, and $Y$ is Moishezon, then the binding group of~$p$ cannot be definably isomorphic to the $\C$-points of a linear algebraic group, unless it is the trivial group.
This is precisely the content of~\cite[Proposition~7.1]{abred}.
Indeed, If~$Y$ is Moishezon and~$f$ is Moishezon, then~$X$ is Moishezon.
But $X$ being Moishezon implies that the $\C$-internality of~$p$ is witnessed without additional parameters, so the binding group is trivial.
\end{remark}

But we can get much more explicit information about types with linear algebraic binding groups, even outside the K\"ahler-type setting:

\begin{proposition}
\label{prop:getvs}
Suppose $K:=\dcl(B)\cap\C$ is algebraically closed and $p\in S(B)$ is a fundamental $\C$-internal and weakly $\C$-orthogonal type whose binding group is definably isomorphic to the $\C$-points of a linear algebraic group.
Then there exists a $B$-definable vector space~$V$ over~$\C$, and a basis~$e$ for~$V$, such that~$p$ is interdefinable with~$\tp(e/B)$.
\end{proposition}

\begin{proof}
By Corollary~\ref{cor:bgag}, we know that there is an algebraic group~$G$ over $K$ such that $G(\C)$ is the dual binding group for~$p$ relative to~$\C$.
As the dual binding group is definably isomorphic to the binding group, our assumption implies that $G$ is a linear algebraic group.
As we are working over the algebraically closed field~$K$, we may assume that  $G\leq\GL_n$ for some $n\geq 1$.
So we write elements of $G(\C)$ as matrices $M$, and denote the dual binding group right action on $P:=p(\overline\M)$, which is a uniquely transitive action, as $aM$ for $a\models p$ and $M\in G(\C)$.
We also have the natural action of $\GL_n(\C)$, and hence of $G(\C)$, on $\C^n$ by linear transformations, which we denote by $Mc$ where $c$ is viewed a a column vector.

Consider the $B$-definable binary relation~$E$ on $P\times\C^n$ given by $(a,c)E(a',c')$ if $Mc'=c$, where $M\in G(\C)$ is the unique element such that $aM=a'$.
It is easily checked that this is an equivalence relation: the identity matrix witnesses reflexivity, the inverse matrix witnesses symmetry.
For transitivity, if $(a_1,c_1)E(a_2,c_2)$ is witnessed by  $a_1M=a_2, Mc_2=c_1$ and $(a_2,c_2)E(a_3,c_3)$ is witnessed by $a_2N=a_3, Nc_3=a_2$, then $a_1MN=a_3$ and $MNc_3=c_1$ so that $(a_1,c_1)E(a_3,c_3)$.

Let $V:=(P\times\C^n)/E$ be the set of $E$-classes, where we denote the $E$-class of $(a,c)$ by $[(a,c)]$.
Note that for any $a,a'\in P$ and $c\in\C^n$ there is a unique $c'\in\C^n$ such that $[(a,c)]=[(a',c')]$.
Indeed, if $M\in G(\C)$ is such that $a'M=a$ then set $c':=Mc$.
Using this fact, we see that addition on $\C^n$ descends to a group structure on $V$ as follows:
Given $[(a,c)],[(b,d)]\in V$ write $[(b,d)]=[(a,d')]$ for some $d'\in\C^n$ and set
$[(a,c)]+[(b,d)]:=[(a,c+d')]$.
We leave it to the reader to verify that this is a well defined additive group structure with identity represented by $(a,0)$ for any $a\in P$.
The induced scalar multiplication is even clearer: $\lambda [(a,c)]:=[(a,\lambda c)]$ is well defined and makes $V$ into a $\C$-vector space.

Let ${\bf e}_1,\dots, {\bf e}_n \in \C^n$ be the standard basis.
Fix $a\in P$.
We leave it to the reader to verify that 
$e_1:=[(a,{\bf e}_1)],\dots, e_n:=[(a,{\bf e}_n)]$
is a $\C$-basis for $V$ (that depends on the choice of~$a$).
It remains to observe that $a$ and $e:=(e_1,\dots,e_n)$ are interdefinable over~$B$.
If $\sigma\in\aut_{Ba}(\overline\M)$ then
\begin{eqnarray*}
\sigma(e_i)
&=&
\sigma([(a,{\bf e}_i)]\\
&=&
[(\sigma(a),\sigma({\bf e}_i))]\ \ \text{ as $E$ is $B$-definable}\\
&=&
[(a,{\bf e}_i]\ \ \text{ as ${\bf e}_i\in\M$ and hence fixed by every automorphism of $\overline\M$}\\
&=&
e_i
\end{eqnarray*}
for all $i\leq n$.
Conversely, if $\sigma\in\aut_{Be}(\overline\M)$ then
$[(a,{\bf e}_i)]=[(\sigma(a),{\bf e}_i)]$ for all~$i$, again because every automorphism fixes each ${\bf e}_i$.
But this means that if $M\in G(\C)$ is such that $aM=\sigma(a)$ then $M{\bf e}_i={\bf e}_i$ for all $i\leq n$.
But this says precisely that $M$ is the identity matrix.
Hence $\sigma(a)=a$, as desired.
\end{proof}

Proposition~\ref{prop:getvs} is in fact a charactisation of types with linear algebraic binding groups.
Indeed, every definable $\C$-vector space has a linear algebraic binding group:

\begin{proposition}
\label{vs-linalg}
Suppose $V$ is a $B$-definable $\C$-vector space, $\Phi$ is a partial type over~$B$ that extends~$V$, and $\G$ is a binding group for~$\Phi$ relative to~$\C$ over~$B$.
Then~$\G$ is definably isomorphic to the $\C$-points of a linear algebraic group.
\end{proposition}
 
 \begin{proof}
As $V$ is of finite rank it must be finite-dimensional as a $\C$-vector space. Fix a $\C$-basis $e=(e_1,\dots,e_n)$ of $V$, and consider $q:=\tp(e/B)$.
Note that~$q$ is (fundamental) $\C$-internal as $q(\overline\M)\subseteq V^n\subseteq\dcl(Be\C)$.
We begin by showing that its binding group is definably isomorphic to the $\C$-points of a linear algebraic group.

First of all, consider the automorphism group $\aut_B(q/\C)$.
Given $\sigma\in \aut_B(q/\C)$, $\sigma(e)$ is another $\C$-basis for $V$.
Let us write $e^t$ for
$\begin{bmatrix}
e_1\\
\vdots\\
e_n
\end{bmatrix}$, and similarly $\sigma(e)^t$.
There is a unique $M_\sigma\in\GL_n(\C)$ such that $M_\sigma e^t=\sigma(e)^t$.
Define $\alpha:\aut_b(q/\C)\to \GL_n(\C)$ by setting $\alpha(\sigma):=M_\sigma^t$, the transpose of $M_\sigma$.
Let us check that $\alpha$ is a group embedding.
Given $\sigma, \tau\in \aut_B(q/\C)$, compute
\begin{eqnarray*}
M_\sigma M_\tau e^t
&=&
M_\sigma\tau(e)^t\\
&=&
M_\sigma
\begin{bmatrix}
\widehat\tau e_1\\
\vdots\\
\widehat\tau e_n
\end{bmatrix}\ \ \ \text{ where $\widehat\tau$ lift $\tau$ to $\aut_{B\C}(\overline\M)$}\\
&=&
\widehat\tau\left(M_\sigma
\begin{bmatrix}
e_1\\
\vdots\\
e_n
\end{bmatrix}
\right)\ \ \text{ as $\widehat\tau$ is the identity on~$B\C$ and $M_\sigma\in\GL_n(\C)$}\\
&=&
\widehat\tau(M_\sigma e^t)\\
&=& \widehat\tau(\sigma(e)^t)\\
&=&(\tau(\sigma(e)))^t\\
&=&
((\tau\sigma)(e))^t\\
&=&M_{\tau\sigma} e^t.
\end{eqnarray*}
As $e$ is a $\C$-basis, it follows that $M_\sigma M_\tau=M_{\tau\sigma}$.
Taking transposes, we get that $\alpha(\tau\sigma)=\alpha(\tau)\alpha(\sigma)$.
Injectivity follows from the fact that $q(\overline\M)\subseteq V^n\subseteq\dcl(be\C)$.
Indeed, if $\alpha(\sigma)=\id$ then $\sigma(e)=e$ and hence $\sigma$ is the identity on all of $q(\overline\M)$.

Let $\H$ be a binding group for~$q$ relative to~$\C$.
Then $\alpha$ induces an embedding $\H\to\GL_n(\C)$.
It is $Be$-definable as one can compute the image of $h\in\H$ by writing the elements of $h(e)$ as $\C$-linear combinations of the $\C$-basis $e$.

Now, let $\G$ be a binding group for~$\Phi$ relative to~$\C$ over~$B$.
As linear algebraic groups are preserved under quotients, it suffices to show that~$\G$ is the image of~$\H$ under a definable homomorphism.
We have surjective homomorphisms from $\aut_{B\C}(\overline\M)$ to $\aut_B(\Phi/\C)$ and to $\aut_B(q/\C)$, given by restriction.
These induce a surjective homomorphism $\aut_B(q/\C)\to\aut_B(\Phi/\C)$.
Indeed, if $\sigma\in\aut_{B\C}(\overline\M)$ acts trivially on $q(\overline\M)$ then it fixes the $\C$-basis~$e$ of $V$, and hence acts trivially on $V$, and so also on $\Phi(\overline\M)$.
We obtain an induced surjective homomorphism $\beta:\H\to\G$.
To see that $\beta$ is definable, note that we have a natural action of  $h\in\H$ on~$V$ given by taking a $\C$-linear combination of $e_1,\dots,e_n$ to the same $\C$-linear combination of the corresponding basis vectors in $h(e)$.
And $\beta(h)=g$ if and only if~$h$ and~$g$ agree on $\Phi(\overline\M)\subseteq V$.
Since the action of~$\G$ on $\Phi(\overline\M)$ is determined by its action on some fixed finite subset of $\Phi(\overline\M)$, namely on a fundamental system of solutions, we have that $\beta$ is definable (over $B$ together with this fundamental system).
\end{proof}

\bigskip
\section{Nontrivial torsors for algebraic groups in $\dccm$}
\label{sect:h1}

\noindent
The theory $\dccm$ was introduced in~\cite{dccm} as a model theoretic context in which to study meromorphic vector fields on meromorphic varieties, in the same way that $\dcf_0$ captures the  model-theoretic approach to rational vector fields on algebraic varieties.
Our point of view is that $\dccm$ expands $\ccm$ in precisely the same way that $\dcf_0$ expands $\acf_0$.
As such it is natural to ask for generalisations of theorems about $\dcf_0$ that refer to it's $\acf_0$-reduct.
One such theorem, that plays a significant role in differential Galois theory, is a theorem of Kolchin's from~\cite{dag}
which says that over an algebraically closed differential field every $\dcf_0$-definable torsor for an algebraic group is trivial.
That is, working in $\dcf_0$, if $G$ is an algebraic group defined over an algebraically closed differential field $K$, and~$X$ is a $K$-definable set admitting a uniquely transitive $K$-definable $G$-action, then~$X$ has a $K$-rational point.
Note that, in general, a $K$-definable set in $\dcf_0$ need not have a $K$-rational point, but this theorem says that if the definable set is a torsor for an {\em algebraic} group, then it does.
See~\cite[Proposition~3.2]{dgt1} for a model-theoretic account of this theorem.
As meromorphic groups\footnote{See~\cite{mergroup} for meromorphic groups in the standard model of $\ccm$, and~\cite{mergroup-nonstandard} for its extension to the theory $\ccm$.}
are the $\ccm$ generalisation of algebraic groups, one might expect the $\dccm$ generalisation of Kolchin's theorem to say something like: over an $\acl$-closed set of parameters every $\dccm$-definable torsor for a meromorphic group is trivial.\
We will show that this is false even for algebraic groups: there are $\dccm$-definable torsors of algebraic groups, over $\acl$-closed parameter sets, that are not trivial.
Indeed, these exist for {\em every} algebraic group that extends a nontrivial $\ccm$-Galois group.

\medskip
\subsection{The $\dccm$ setting}
Let $L$ be, as before, the language of meromorphic varieties, and set $L_\nabla$ to be the expansion of $L$ obtained by adding function symbols $\nabla_S:S\to TS$ for each sort~$S$.
Here $TS\to S$ is the tangent space of the meromorphic variety~$S$.
Then $\dccm$, the $L_\nabla$-theory of {\em differentially closed $\ccm$-structures}, axiomatises the existentially closed ``differential $\ccm$-structures".
We leave the reader to consult~\cite{dccm} for more precise definitions and details.
In particular, it is shown there that $\dccm$ is a complete totally transcendental theory admitting the eliminations of quantifiers and imaginaries.

We work in a sufficiently saturated model $(\U,\nabla)\models\dccm$.

Let $\overline\M=\mathcal U^\nabla$ denote the {\em constants} of $(\U,\nabla)$, namely the $L$-structure where for each sort~$S$ we set $S(\overline\M):=\{a\in S(\U):\nabla(a)=0\}$.
Then $\overline\M$ is a saturated model of $\ccm$ that is definable in $(\U,\nabla)$, and the $L$-structure agrees with full induced structure from $(\U,\nabla)$.
In this way, $\ccm$ lives as a pure stably embedded reduct of $\dccm$.

Let $\K:=\AA(\U)$, it is an algebraically closed field extending~$\CC$.
Observe that $\nabla_\AA$ induces a derivation~$\delta$ on $\K$ such that $(\K,\delta)$ is a saturated model of $\dcf_0$ definable in $(\U,\nabla)$.
Moreover, the full induced structure on $\K$ from $(\U,\nabla)$ is interdefinable with $(\K,\delta)$.
In this way $\dcf_0$ lives as a pure stably embedded reduct of $\dccm$.

The field of constants of $(K,\delta)$ is $\C:=\AA(\overline\M)$.

We thus have three models of $\ccm$ here:
$\dcl(\emptyset)=\M\subseteq\overline\M\subseteq\U$.
And if we take their $\AA$-sorts we get three algebraically closed fields: $\CC\subseteq\C\subseteq\K$.

\begin{remark}
We introduce some conventions and terminology:
\begin{itemize}
\item[(i)]
We view $\dcl$-closed sets, $B$, as many-sorted substructures; so for each sort~$S$ we have $S(B)\subseteq S(\U)$.
In particular, $\K_B:=\AA(B)$ is a $\delta$-subfield of $\K$, and $\C_B:=\K_B\cap\C$ is its field of constants.
\item[(ii)]
We identify any (quasi-projective) variety~$V$ over~$\K$ with its $\K$-rational points, $V(\K)$, which is itself identified with a definable set in the field~$\K$.
In particular, varieties are thus viewed as definable objects in $(\U,\nabla)$.
\item[(iii)]
It follows from~(i) and (ii), that if $V$ is a variety over~$\K$, and~$B$ is a $\dcl$-closed set, then it makes sense to consider the set of $B$-points of $V$, denoted by $V(B)$, and that this agrees precisely with the set of $\K_B$-rational points of $V$.
Moreover, because~$\K$ is stably embedded, $V$ will be defined over~$B$ (in the model-theoretric sense) if and only if it is defined over the $\delta$-field $\K_B$ (in both the model-theoretic and algebro-geometric senses).
\end{itemize}
\end{remark}

\medskip
\subsection{Existence of a nontrivial torsor}

\begin{theorem}
\label{thm:h1}
Suppose $G$ is an infinite $\ccm$-Galois group over~$\CC$.
Then $G$ admits a nontrivial definable torsor in $\dccm$, over an $\acl$-closed set of parameters.

Moreover, the same is true of any algebraic group into which $G$ embeds.
\end{theorem}

\begin{proof}
Our proof makes use of~\cite{galco}, which describes definable torsors in terms of definable cocycles.
So what we construct is a nontrivial definable cocycle.

We first work in $\overline\M\models\ccm$, with algebraically closed field $\C=\AA(\overline\M)$ extending~$\CC$.
That $G$ is a $\ccm$-Galois group means, using Lemma~\ref{galois-dual}, that there is $B=\acl_L(B)$, and $p\in S(B)$ fundamental $\C$-internal and weakly $\C$-orthogonal, such that $G(\C)$ is the dual binding group of~$p$ relative to $\C$.

Now pass to the expansion $(\U,\nabla)\models\dccm$.
So $\overline\M=\U^\nabla$.
Because the full induced structure on~$\overline\M$ by $(\U,\nabla)$ is just the $L$-structure,
we have that $p(x)$ together with the formula $\nabla(x)=0$ is a complete type over~$B$ in $(\U,\nabla)$, that we still call~$p$, and such that $p(\U)=p(\overline\M)$.
Moreover, $p$ remains fundamental $\C$-internal and weakly $\C$-orthogonal in $(\U,\nabla)$.
And $G(\C)$ remains the dual binding group of~$p$ relative to~$\C$ when computed in $(\U,\nabla)$.
As $\overline\M$ is $\acl$-closed in~$(\U,\nabla)$, so is~$B$.

Let $\widehat B$ be a prime model of $\dccm$ over~$B$, and fix $a\in p(\widehat B)$.
This exists as $p$ is isolated.
Let $D:=\dcl(Ba)\subseteq\widehat B$.

Note that $\C_{\widehat B}=\C_B$ (so that we can apply Lemma~\ref{fact:widehatbtod}).
This is because $B$ is $\acl$-closed, and hence $\C_B$ is an algebraically closed field, but $\C_{\widehat B}$ is algebraic over $\C_B$ as $\widehat B$ is prime over~$B$.

Composing the homomorphism from Corollary~\ref{cor:extrinsic}(b) and Lemma~\ref{fact:widehatbtod}, applied to $p$ and its dual binding group $G(\C)$ in $\dccm$, gives a surjective homomorphism
$$\aut_B(\widehat B)\to G(B).$$
Since $G(B)\leq G(\widehat B)$, this is a homomorphism to the latter.
That is, we have a homomorphism
$$s:\aut_B(\widehat B)\to G(\widehat B),$$
whose image lands in $G(B)$, and which 
satisfies
$as(\sigma)=\sigma(a)$
for all $\sigma\in\aut_B(\widehat B)$.

We claim that $s$ is a {\em cocycle}: namely that $s(\sigma_1\sigma_2)=s(\sigma_1)\sigma_1(s(\sigma_2))$ for all $\sigma_1,\sigma_2\in \aut_B(\widehat B)$.
Indeed, this follos from the fact that $s$ is a homomorphism, and its image lands in $G(B)$ which is fixed pointwise by $\aut_B(\widehat B)$.

Let $h(x,y)$ be the $B$-definable function on $p(x)\times p(y)$ that assigns to $(x,y)$ the unique element of $G(\C)$ that takes $x$ to $y$.
Then $s(\sigma)=h(a,\sigma(a))$, for all $\sigma\in\aut_B(\widehat B)$.
That is, $s$ is a {\em definable cocycle} in the sense of~\cite{galco}.
(Note that~$h$ is in fact $L$-definable as its graph lives in the stably embedded $\overline\M$, but we don't make use of this fact.)

Now we apply the main result of~\cite{galco}, which gives a natural correspondence between the set of definable cocycles, modulo cohomology, and the $B$-definable right torsors for $G$, modulo isomorphism.
In particular, the torsor associated to $s$ will be trivial (that is, have a $B$-point) if and only if $s$ is the {\em trivial cocycle} in the sense that there is $\alpha\in G(\widehat B)$ such that $s(\sigma)=\alpha^{-1}\sigma(\alpha)$, for all $\sigma\in\aut_B(\widehat B)$.
So it remains to show that $s$ is nontrivial as a cocyle

Suppose, toward a contradiction, that there is $\alpha\in G(\widehat B)$ with $s(\sigma)=\alpha^{-1}\sigma(\alpha)$ for all $\sigma\in\aut_B(\widehat B)$.

\begin{claim}
\label{indclBa}
$\alpha\in D$.
\end{claim}

\begin{claimproof}
Suppose $\alpha'\models\tp(\alpha/Ba)$ is in $\widehat B$.
Then there is $\sigma\in\aut_{Ba}(\widehat B)$ such that $\sigma(\alpha)=\alpha'$.
But
$$1_G=h(a,a)=h(a,\sigma(a))=s(\sigma)=\alpha^{-1}\sigma(\alpha)=\alpha^{-1}\alpha'.$$
That is, $\alpha'=\alpha$.
So $\tp(\alpha/Ba)$ has only one realisation in $\widehat B$, and as the type is isolated (as we are in the prime model over $B$), we have $\alpha\in\dcl(Ba)=D$.
\end{claimproof}

\begin{claim}
\label{incb}
$G(D)=G(B)$
\end{claim}

\begin{claimproof}
Since $B\subseteq\overline\M=\U^\nabla$, and $\nabla(a)=0$, we have that
$$D=\dcl(Ba)\subseteq\overline\M$$
as well.
On the other hand, $G$ lives on the sort~$\AA$, so that $G(\overline\M)=G(\C)$.
Hence $G(D)= G(\C_D)$.
Finally, as $\tp(a/B)$ is weakly orthogonal to~$\C$ in $\overline\M$, we have that $\C_D=\C_B$ by Lemma~\ref{lem:wochar}.
\end{claimproof}

So $\alpha\in G(B)$, and we have that $\sigma(\alpha)=\alpha$ for all $\sigma\in\aut_B(\widehat B)$.
It follows that $s(\sigma)=\alpha^{-1}\sigma(\alpha)=1_G$ is the trivial homomorphism.
Hence its image, $G(B)$, is the trivial group.
As $\C_B$ is an algebraically closed field over which~$G$ is defined, this contradicts the positive-dimensionality of~$G$.

The ``moreover" clause of Theorem~\ref{thm:h1} follows from an inspection of the above proof: If $G'$ is an algebraic group in which $G$ embeds (over~$\CC$), then we can view $s$ as a homomorphism $\aut_B(\widehat B)\to G'(\widehat B)$ as well.
The proof that $s$ is a nontrivial definable cocycle simply goes through with $G'$ in place of $G$.
And we thus obtain a nontrivial $B$-definable torsor for $G'$ in $\dccm$, as desired.
\end{proof}

One case of Theorem~\ref{thm:h1} is that every $\GL_n$ has a nontrivial definable torsor in $\dccm$.
Indeed, we have seen that $\mathbb G_m$ is a $\ccm$-Galois group, and $\mathbb G_m$ embeds in each $\GL_{n}$.
This can be seen as saying that Hilbert's 90th does not hold in $\dccm$.

\medskip
\subsection{An explicit nontrivial torsor}
We now trace through the correspondence in~\cite{galco}, between cocycles and right torsors, to describe explicitly a nontrivial right torsor for any given algebraic group $G$ arising as a $\ccm$-Galois group.

We work with the construction in the proof of Theorem~\ref{thm:h1}.
In particular, $p$ is the fundamental $\C$-internal and weakly $\C$-orthogonal type in $\dccm$, with dual binding group $G(\C)$, from that proof.
Let $P$ denote the set of realisations of $p$, it is a $B$-definable set by isolation of~$p$.
(Note that $P\subseteq\overline\M$.)
Let $h:P^2\to G(\C)$ be the $B$-definable function from that construction, so $h(x,y)$ is the unique element of $G(\C)$ that sends $x$ to $y$.

\begin{claim}
\label{2.7}
For all $a,b,c\in P$, $h(a,c)=h(a,b)h(b,c)$.
\end{claim}

\begin{claimproof}
It suffices to check this in the prime model $\widehat B$ over $B$.
But there, we have that $\sigma\mapsto h(a,\sigma(a))$ defines a cocycle from $\aut_B(\widehat B)\to G(\widehat B)$.
The above identity is then precisely Lemma~2.7 of~\cite{galco}, and is easily verified using the fact that any two elements of $P(\widehat B)$ are $\aut_B(\widehat B)$-conjugate.
\end{claimproof}

Now consider the $B$-definable set $P\times G$, and let $E$ be the $B$-definable relation on $P\times G$ given by
$$(c,\alpha)E(c',\alpha')\iff\alpha=h(c,c')\alpha'.$$
Using Claim~\ref{2.7}, and fixing $a\in P$, we see that 
$$(c,\alpha)E(c',\alpha')\iff h(a,c)\alpha=h(a,c')\alpha'.$$
It follows immediately that~$E$ is an equivalence relation on $P\times G$.

Consider the $B$-definable set $X:=(P\times G)/E$.
It is clear that the right action of~$G$ on $P\times G$, given by  $(c,\alpha)\cdot g:=(c,\alpha g)$, descends to a $B$-definable action on $X$.
The action is transitive: given $(c,\alpha),(c',\alpha')\in P\times G$, if we set $g=\alpha^{-1}h(c,c')\alpha'$ then $\alpha g=h(c,c')\alpha'$, so that $(c,\alpha)\cdot g$ is $E$-equivalent to  $(c',\alpha')$.
It is uniquely transitive: if $(c,\alpha)\cdot g$ is $E$-equivalent to $(c,\alpha)$ then $\alpha g=h(c,c)\alpha=\alpha$ (as, by Claim~\ref{2.7}, for example, $h(c,c)=1_G$), and hence $g=1_G$.
So $X$ is a $B$-definable right torsor for~$G$.
The nontriviality of the cocycle says precisely that $X$ has no $B$-points.

\vfill\pagebreak

\appendix

\section{Complex tori with prescribed algebraic reduction\\ by Fr\'ed\'eric Campana}
\label{campana-appendix}

\noindent
At the end of Section~4 of~\cite{cored} it is stated that complex tori can be constructed with various prescribed properties and that proofs will be given later.
Here we prove a very special case of this claim, following a question of the authors of this paper.
Remember that an abelian variety is precisely a complex torus whose underlying variety is projective.

\begin{proposition}
\label{campana}
Let $A, B$ be complex abelian varieties with $A$ simple.
Then there is a complex torus~$T$ and a short exact sequence 
$$
\xymatrix{
0\ar[r] & A\ar[r] &T\ar[r]^p & B\ar[r] & 0
}$$
such that $p:T\to B$ is the algebraic reduction of $T$. 
\end{proposition}

We begin with a general discussion of arbitrary short exact sequences
$$
\xymatrix{
0\ar[r] & A\ar[r] &T\ar[r]^p & B\ar[r] & 0
}$$
of complex tori.
See~\cite{tori} for details.
Let $d = \dim(B)$. 
Note that $p:T\to B$ is a principal fibre bundle.
Let $\pi_{1}(B)$ be the fundamental group of~$B$.
So $\pi_{1}(B) = \mathbb Z^{2d}$ with an embedding in the universal cover
$\widetilde B= \CC^{d}$ of $B$, with quotient~$B$.
On the other hand we have the monodromy representation~$r$ of~$\pi_{1}(B)$ in $\aut^{0}(A)= A$ (acting by translations).
To remind the reader, $r$ is defined as follows:
Take a path~$\gamma$ in~$B$ starting and ending in the identity $0=0_{B}$.
Let $x\in A = p^{-1}(0)$.
Then~$\gamma$ lifts to a path in~$T$ from~$x$ to some~$y\in A$ and we write $y = r(\gamma)(x)$.
One can check directly that~$r:\pi_1(B)\to A$ is a homomorphism (of abstract groups). 
And that~$T$ can be recovered as the quotient of
${\tilde B}\times A$ by the discrete subgroup $\Gamma:=\{(\gamma, r(\gamma)):\gamma\in \pi_{1}(B)\}$. 

Conversely, starting with complex tori $A$ and $B$, and any (abstract) homomorphism $r:\pi_{1}(B) \to A$, we obtain a complex torus extending~$B$ by~$A$  with monodromy representation~$r$.
Indeed, set $T:=({\tilde B}\times A)/\Gamma$ where $\Gamma$ is as above.

\begin{lemma}
\label{av=fin}
Let
$\xymatrix{
0\ar[r] & A\ar[r] &T\ar[r]^p & B\ar[r] & 0}$ be a short exact sequence of complex tori, where $A,B$ are abelian varieties, and let $r:\pi_{1}(B)\to A$ be the associated monodromy representation.
Then
$T$ is an abelian variety if and only if the image of~$r$ is finite.
\end{lemma}

\begin{proof} If $T$ is an abelian variety then by Poincar\'e complete reducibility there is some abelian subvariety $B'$ of $T$ such that $p$ induces an isogeny $B' \to B$, and in particular $A\cap B'$ is finite.
One checks that the image of $r$ is in $A\cap B'$.
For the converse,
if the image of~$r$ in~$A$ is finite then  $T=({\widetilde B}\times A) /\Gamma$ is isogeneous to $B\times A$, and hence projective.
\end{proof}

Hence by the discussion preceding the lemma, starting with complex abelian varieties $A, B$ and a homomorphism $r:\pi_{1}(B)\to A$ with infinite image  (for example some generator of $\pi_{1}(B)$ is sent to a torsion-free element of $A$), we can find, by the lemma, a nonalgebraic complex torus~$T$ with $A\leq T$ and $T/A = B$.
As the algebraic reduction $f: T\to G$ of~$T$ must factor through~$p:T\to B$, 
we have that $\ker(f)\leq A$. 
So if~$A$ is a simple abelian variety, then $\ker(f)$ is either trivial or all of~$A$.
It cannot be trivial as~$T$ is nonalgebraic.
Hence $G = B$ and $p:T\to B$ is the algebraic reduction.
This proves Proposition~\ref{campana}.
\qed

For completeness we present a generalisation of Lemma~\ref{av=fin} that drops the assumption that~$A$ and~$B$ are abelian varieties.

\begin{lemma}
Let
$\xymatrix{
0\ar[r] & A\ar[r] &T\ar[r]^p & B\ar[r] & 0}$ be a short exact sequence of complex tori.
Then the following are equivalent:
\begin{itemize}
\item[(i)]
There exists a Zariski closed subset $X\subseteq T$ such that $p|_X: X \to B$ is surjective and generically finite-to-one.
\item[(ii)]
There is a holomorphic map $f:B \to T$ such that $p\circ f$ is an isogeny.
\item[(iii)]
There is an isogeny $g: B\times A \to T$ such that $p\circ g$ restricts to an isogeny $B\times \{0_{A}\} \to B$. 
\item [(iv)]
The associated monodromy representation $r:\pi_{1}(B)\to A$ has finite image.
\end{itemize}
\end{lemma}

\begin{proof} We give a sketch.

(i)$\implies$(ii).
For generic $u\in B$, define $f(u)  = \sum\{x\in X: p(x) = u\}$, with sum  taken in the group $T$.
This defines a meromorphic map $f: B\to T$ which must therefore be holomorphic.
So $p\circ f:B\to B$ is holomorphic and agrees with multiplication by~$m$ generically (where $m$ is the size of the generic fibre of~$X$ over ~$B$), and hence everywhere.

(ii)$\implies$(iii).
Define $g(b,a) = f(b) + a$.

(iii)$\implies$(iv).
One checks that image of $r$ is in $A\cap B'$ where $B':= g(B\times\{0_A\})$.

(iv)$\implies$(i).
As in the proof of Lemma~\ref{av=fin}, we get that $T$ is isogenous to $B\times A$ over~$B$ and we can take $X$ to be the image of $B\times\{0_A\}$.
\end{proof}

\vfill\pagebreak

\section{Strong irreducibility of the tangent bundle of an\\irreducible hyperkaehler manifold\\ by Matei Toma}
\label{toma-appendix}

\noindent
We will consider here irreducible compact hyperk\"ahler manifolds $Y$ of (complex) dimension $2n$ and their projectivized tangent bundles $\PP(T_{Y})\to Y$. We will show that under some condition on $Y$, which can be met for all $n>0$, there exist no proper Zariski closed subspaces in $\PP(T_{Y})$ which project onto $Y$. Here $\PP(T_{Y})\to Y$ denotes the $\PP^{2n-1}$-bundle of hyperplanes in the fibers of $T_{Y}$ over $Y$.  Yau's solution to the Calabi Conjecture and a classical result of Weyl on the representations of the compact symplectic group will play an important role in the proof.

A compact complex manifold $Y$ of dimension $2n$ is said to be {\em irreducible hyperk\"ahler} if it admits a K\"ahler metric whose holonomy equals the compact symplectic group $\Sp(n)$.  This condition is equivalent, by Yau's solution to the Calabi Conjecture, to asking that $Y$ be an {\em irreducible holomorphic symplectic manifold} by which is meant that $Y$ is K\"ahler, simply-connected and admits an everywhere non-degenerate holomorphic two-form, which is unique up to a constant factor, \cite[Proposition 4]{Beauville1983}. See for instance \cite{Beauville1983, GrossHuybrechtsJoyce} for their description and properties.

\begin{proposition}
\begin{enumerate}
\item[(a)] If  $Y$ is an irreducible compact hyperk\"ahler manifold then all the symmetric powers  $S^m T_Y$ of its tangent bundles are stable with respect to any K\"ahler class $[\omega]$ on $Y$. This means that for any non-trivial proper coherent subsheaf $\F$ of $S^m T_Y$ with $\rank(\F)<\rank(S^m T_Y)$ one has $\deg_{[\omega]}(\F):=\int_{X}c_{1}(\F)\wedge[\omega]^{2n-1}<0.$
\item[(b)] If in the above situation we suppose in addition that the Picard group $\Pic(Y)$ is trivial, 
then the projectivized tangent bundle 
 $\PP(T_Y)$ of $Y$ admits no proper Zariski closed subspaces projecting onto $Y$. In particular in this case both $Y$ and $\PP(T_Y)$ have zero algebraic dimension. 
\item[(c)] For all $n>0$ there exist irreducible compact hyperk\"ahler manifolds $Y$ of dimension $2n$ with trivial Picard group.
\end{enumerate}
\end{proposition}

\begin{proof}
(a).
Let $Y$ be an irreducible compact hyperk\"ahler manifold  of dimension $2n$ and let $[\omega]$ be a K\"ahler class on $Y$. In this case Yau's solution to the Calabi Conjecture implies that there exists a representative of this class (which we will denote again by $\omega$) such that  the holonomy of the associated Riemannian structure on $Y$ is exactly $\Sp(n)$.  It also implies that $(Y,\omega)$ is K\"ahler-Einstein and that  the induced Hermitian metric on the tangent bundle  $T_Y$ is irreducible Hermite-Einstein with respect to $\omega$.  The irreducibility is due to the fact that the standard $\Sp(n)$ representation on $\CC^{2n}$ is irreducible.  Now by a classical result from \cite{Weyl}  all symmetric
powers $S^m(\CC^{2n})$  of the standard representation are  irreducible,  cf.  \cite[§24.2,  p. 406]{FultonHarris}.  
This implies that the induced Hermitian metrics on the symmetric powers 
 $S^mT_Y$ of the tangent bundle of $Y$ are irreducible Hermite-Einstein with respect to~$\omega$.  By the Kobayashi-Hitchin correspondence it follows that the bundles  $S^mT_Y$ are all $[\omega]$-stable (and not just $[\omega]$-polystable).

(b).
We start by an easy observation.
If $X$ is a compact complex manifold with trivial Picard group and $E$ is a torsion-free stable coherent sheaf on $X$ (with respect to some Gauduchon metric), then $E$ is {\em irreducible}, i.e. it admits no coherent subsheaf $\F$ of intermediate rank. 
Indeed, the triviality of the Picard group implies that all torsion-free coherent sheaves have degree zero with respect to any Gauduchon metric on $X$, so stability immediately implies irreducibility.

Thus by (a) for an  irreducible compact hyperk\"ahler manifold $Y$ with vanishing Picard group, all the symmetric powers $S^mT_Y$  of the tangent bundle are irreducible holomorphic vector bundles.  
Consider now in this situation the projection  $\pi:\PP(T_Y)\to Y$,  and suppose by contradiction that a proper closed subvariety  $Z$ which covers $Y$ existed in $\PP(T_Y)$. Denote by $\O_{ \PP(T_Y)}(1)$ the tautological line bundle on $\PP(T_Y)$, by $F$ the general fibre of  $Z\to Y$ and
by $\I_F$ its ideal sheaf inside the fibre $\PP^{2n-1}$ of $\pi$ which contains $F$. We have an inclusion of coherent sheaves $\I_F \subset \O_{\PP^{2n-1}}$. 
Tensoring it by  the tensor power $\O_{\PP^{2n-1}}(1)^{\otimes m}$ of the tautological line bundle on $\PP^{2n-1}$ for $m\ge0$ and taking global sections we get
$$H^0(\PP^{2n-1},\I_F(m))\varsubsetneq H^0(\PP^{2n-1},\O_{\PP^{2n-1}}(m)).$$
Moreover for $m$ large enough one has $H^0(\PP^{2n-1},\I_F(m))\ne0$ which gives rise to an inclusion of coherent sheaves 
$$0\ne\pi_* (\I_Z(m)) \varsubsetneq \pi_* ( \O_{ \PP(T_Y)}(m))=S^mT_Y$$ on $Y$
such that $0<\rank (\pi_*(\I_Z(m))<\rank( S^mT_Y)$ and contradicts the irreduciblity of $S^mT_Y$.

The triviality of the Picard group of $Y$ implies that $Y$ admits no closed hypersurfaces and therefore that $a(Y)=0$.  Further, $\PP(T_Y)$ cannot admit closed hypersurfaces either, since by dimension reasons they would cover either $Y$ or some hypersurface of $Y$. So the algebraic dimension of $\PP(T_Y)$ vanishes too.

(c). By  \cite[Section 8 Remarque 1] {Beauville1983} 
small deformations of an  irreducible compact hyperk\"ahler manifold remain  irreducible hyperk\"ahler.  Moreover, using the period map, the local Torelli Theorem, the  Beauville-Bogomolov-Fujiki intersection form and \cite[Section 8 Corollaire 1] {Beauville1983} 
one can show that there exist arbitrarily small deformations of any given irreducible compact hyperk\"ahler manifold which will have trivial Picard group.  

The fact that there exist irreducible compact hyperk\"ahler manifolds of any even dimension $2n$ for $n>1$ is established by Beauville's series of examples \cite{Beauville1983}. 
\end{proof}

%\bibliographystyle{plain}
%\bibliography{h1}

\end{document}